\documentclass[11pt, leqno]{amsart} 
\usepackage{amsthm, amsmath}
\usepackage{amssymb}
\usepackage{amscd}
\usepackage[all]{xy}
\usepackage{verbatim}
\usepackage{pdfsync}
\usepackage{graphicx}
\usepackage{color}
\usepackage{tikz}

\usepackage{marginnote}

\newcounter{mynote}

\theoremstyle{plain}
\newtheorem{thm}{Theorem}[section]

\newtheorem{lemma}[thm]{Lemma}

\newtheorem{cor}[thm]{Corollary}

\newtheorem{Question}[thm]{Question}

\theoremstyle{definition}
\newtheorem{dfn}[thm]{Definition}
\newtheorem{example}[thm]{Example}

\theoremstyle{remark}
\newtheorem{rem}[thm]{Remark}

\newcounter{ex}[section]

\newcommand{\Char}{{\rm char}}

\newcommand{\Gal}{{\rm Gal}}

\DeclareMathOperator{\Ker}{Ker}
\DeclareMathOperator{\Hom}{Hom}

\renewcommand{\to}{\longrightarrow}

\newcommand{\mR}{{\mathbb R}}
\newcommand{\mZ}{{\mathbb Z}}

\newcommand{\mQ}{{\mathbb Q}}

\DeclareMathOperator{\Alb}{Alb}
\DeclareMathOperator{\NS}{NS}
\DeclareMathOperator{\Num}{Num}

\newcommand{\piet}{\pi_1^{\textrm{\'et}}}

\begin{document}
\title[Negative curves of small genus on surfaces]{Negative curves of small genus on surfaces}

\author[T. Chinburg]{Ted Chinburg*}\thanks{*Supported in part by NSF Grants DMS 1360767, DMS 1265290 and DMS 1100355, SaTC grant CNS-1513671 and Simons fellowship 338379. The author acknowledges support from U.S. National Science Foundation grants DMS 1107452, 1107263, 1107367 ``RNMS: GEometric structures And Representation varieties" (the GEAR Network).}
\address{Ted Chinburg\\ Department of Mathematics\\University of Pennsylvania\\ Philadelphia, PA 19104, U.S.A.}
\email{ted@math.upenn.edu}
 
\author[M. Stover]{Matthew Stover**}\thanks{**Supported in part by NSF RTG grant DMS 0602191, NSF grant DMS 1361000, and Grant Number 523197 from the Simons Foundation/SFARI. The author acknowledges support from U.S. National Science Foundation grants DMS 1107452, 1107263, 1107367 ``RNMS: GEometric structures And Representation varieties" (the GEAR Network).}
\address{Matthew Stover\\ Department of Mathematics\\Temple University\\Philadelphia, PA 19122, U.S.A.}
\email{mstover@temple.edu}

\date{\today}

\bigskip


\begin{abstract}
Let $X$ be a smooth geometrically irreducible projective surface over a field. In this paper we give an effective upper bound in terms of the N\'eron--Severi rank of $X$ for the number of irreducible curves $C$ on $X$ with negative self-intersection and geometric genus less than $b_1(X)/4$, where $b_1(X)$ is the first \'etale Betti number of $X$. The proof involves a hyperbolic analog of the theory of spherical codes. More specifically, we relate these curves to the \emph{hyperbolic kissing number}, and then prove upper and lower bounds for the hyperbolic kissing number in terms of the classical Euclidean kissing number.
\end{abstract}

\maketitle

\section{Introduction}
\label{sec:intro}

Let $X$ be a smooth geometrically irreducible projective surface over a field $k$. By a negative curve on $X$ we will mean a complete, reduced, irreducible curve with negative self intersection. The bounded negativity conjecture states that if $k = \mathbb{C}$, there is a lower bound depending only on $X$ for the self-intersection of all negative curves on $X$. The origins of this conjecture, as well as various results concerning it, are discussed in \cite{AllThem}, \cite{MollerToledo}, \cite{KoziarzMaubon} and their references. In this paper we study a different but related question that was asked explicitly in the introduction of \cite{AllThem}:
\begin{Question}
\label{eq:basicq}
Let $X$ be a smooth geometrically irreducible projective surface over an arbitrary field $k$. For which integers $g \ge 0$ are there infinitely many negative curves on $X$ of genus $g$?
\end{Question}
\noindent Bogomolov proved in \cite{Bogomolov} that if $k = \mathbb{C}$, $X$ has general type and $c_1(X)^2 > c_2(X)$, the number of negative curves $C$ on $X$ of geometric genus $g$ is finite for every $g$. He did this by showing that, under these hypotheses, the $C$ on $X$ that have a given geometric genus form a bounded family. An effective version of the latter result was shown by Lu and Miyaoka in \cite[Thm.\ 1(1)]{LuMiyaoka}. We refer the reader to \cite[\S 2]{AllThem} for a more detailed account of subsequent related work. The results of Bogomolov and Lu--Miyaoka need not hold for arbitrary $X$ even over $k = \mathbb{C}$. For example, it was shown in \cite[Thm.\ 4.3]{AllThem} that for each $m > 1$ and each $g \ge 0$ there is an $X$ over $k = \mathbb{C}$ containing infinitely many smooth irreducible curves of self-intersection $-m$ and genus $g$.

In this paper we prove that for all $k$ and $X$, there is an effective finite upper bound on the number of negative curves $C$ on $X$ of geometric genus less than $b_1(X)/4$, where $b_1(X)$ is the first \'etale Betti number of $X$. The proof uses a different method than Bogomolov's and involves studying a hyperbolic analog of the theory of spherical codes in Euclidean space.
 
The theory of spherical codes arose from the classical question, going back to Newton and Gregory, of determining how many unit spheres in Euclidean $n$-space can touch a unit sphere centered at the origin without overlapping interiors. The centers of these spheres then form a spherical code on a sphere of radius $2$ about the origin. The angle of such a code is at least $\pi/3$ in the sense that rays from the origin to two such centers must form an angle of at least $\pi/3$. There is an extensive literature concerning spherical codes of varying angles; see, for example, \cite{ConwaySloane} and \cite{Ericson}.

We will show how a hyperbolic analog of this theory arises naturally from studying negative curves of genus less than $b_1(X)/4$ on surfaces. The relevant angle is then $\pi/2$ rather than $\pi/3$, as we will explain below. The study and classification of which hyperbolic codes can in fact arise from negative curves in this way is a natural one when trying to determine the constraints on such curves coming from intersection theory.
 
To state an explicit result, recall that $b_1(X) = \mathrm{dim}_{\mQ_p}\, H^1_{\textrm{\'et}}(X, \mQ_p)$ for any prime $p$ different from $\Char(k)$. The Picard number $\rho(X)$ is the rank over $\mZ$ of the N\'eron--Severi group $\NS(X)$ of $X$. We prove:

\begin{thm}
\label{thm:first}
Let $X$ be a smooth geometrically irreducible projective surface over a field $k$. There are effective constants $c_1$ and $c_2$ for which the number $\tau(X)$ of negative curves $C$ on $X$ with geometric genus $g(C) < b_1(X) / 4$ is finite and bounded above by $c_1 e^{c_2 \rho(X)}$. One has $\tau(X) \le 2^{\, 0.902\, \rho(X)}$ for sufficiently large $\rho(X)$.
\end{thm}

We do not know if the condition $g(C) < b_1(X) / 4$ in Theorem \ref{thm:first} is sharp in characteristic zero. However, it is sharp in positive characteristic in view of the following example.

\begin{example}\label{ex:hell}
Let $X$ be the direct product $Y \times Y$, where $Y$ is a smooth projective irreducible curve of genus $g \ge 2$ defined over $\mathbb{F}_p$. Set $k = \overline{\mathbb{F}}_p$, and let $C_n$ be the graph of $\sigma^n$ in $X(\overline{\mathbb{F}}_p)$, where $\sigma: Y \to Y$ is the Frobenius automorphism. Then $C_n$ is a reduced irreducible curve on $X$ of arithmetic genus $g$, and \cite[Ex.\ V.1.10]{Hartshorne} implies that $C_n^2 \to -\infty$ as $n \to \infty$. One has $b_1(X) = 4 g$. In particular, there are infinitely many distinct reduced irreducible curves on $X$ of genus $g = b_1(X) / 4$ with negative self-intersection.
\end{example}

The proof of Theorem \ref{thm:first} is effective and involves analyzing how negative curves of small genus on $X$ are constrained by the intersection pairing on the Neron--Severi group via the Hodge index theorem. We will precisely define the set of constraints to be considered. In the course of proving Theorem \ref{thm:first}, we will also show that one cannot use these constraints alone to improve the upper bound on $\tau(X)$ to one that is subexponential in $\rho(X)$.

We discuss below how the curves of Theorem \ref{thm:first} give rise to a hyperbolic code of angle at least $\pi/2$. Moreover, we will show that in fact, there do exist hyperbolic codes of angle at least $\pi/2$ that grow exponentially in size with $\rho(X)$. We state this in terms of the \emph{hyperbolic kissing number}. Recall that the classical kissing number $K_n$ arises from finding maximal spherical codes in the unit sphere $\mathbb{S}^{n-1}$ in Euclidean $n$-space which have a given angle. See \S \ref{s:spherical} for more on the classical case. In \S \ref{s:hyperbolic}, we define the hyperbolic kissing number in an analogous manner, and one of our key technical results is Theorem \ref{thm:hyperbounds}, which gives upper and lower bounds for the hyperbolic kissing number in terms of classical kissing numbers.

However, we do not know whether such large codes can arise from the intersection theory of negative curves of small genus on surfaces. This leads to the following question:

\begin{Question}
Is there a constant $c > 1$ and a sequence of surfaces $X$ for which $\rho(X)$ tends to infinity such that $\tau(X) > c^{ \rho(X)}$ for all $X$ in this sequence? Can one find such a sequence in which $b_1(X)$ remains bounded?
\end{Question}

To describe the connection between negative curves and hyperbolic codes more precisely, recall that $\NS(X)$ is a finitely generated abelian group, and the group $\Num(X)$ of divisors modulo numerical equivalence on $X$ is the quotient of $\NS(X)$ by its torsion subgroup. Thus $\Num(X)$ is a free $\mathbb{Z}$-module of rank $\rho(X) \ge 1$. The curves $C$ in Theorem \ref{thm:first} map bijectively to their classes in $\Num(X)_{\mathbb{R}} = \mathbb{R} \otimes_{\mathbb{Z}} \Num(X)$. We show that these classes form a strict hyperbolic code of angle at least $\pi / 2$ in $\Num(X)_{\mathbb{R}}$ in the sense of Definition \ref{def:hypercodedef}. The maximal number of elements in such a code is the strict hyperbolic kissing number $R_{\rho(X) - 1}(\pi / 2)$. We will prove:

\begin{thm}\label{thm:IntroKiss} Let $R_n(\pi/2)$ be the strict hyperbolic kissing number of hyperbolic $n$-space. Then $R_0(\pi/2) = 0$ and $R_1(\pi/2) = 1$. For $n \ge 2$, $R_n(\pi/2)$ is bounded from above by $K_{n-1}(\phi_0) + 2$, where $K_{n-1}(\phi_0)$ is the classical kissing number for the Euclidean unit sphere $\mathbb{S}^{n-2}$ in $\mR^{n-1}$ associated with the angle $\phi_0 = \mathrm{arccos}(3/4)$. Furthermore, if $n \ge 2$ then $R_n(\pi/2)$ is bounded below by $\lfloor K_{n-1}(2\phi_0)/2 \rfloor$, the greatest integer less than or equal to $K_{n-1}(2\phi_0)/2$.
\end{thm}

The bound in Theorem \ref{thm:first} then follows from an upper bound of Kabatiansky and Levenshtein on $K_{n-1}(\phi_0)$. This lower bound grows exponentially with $n$ by work of Chabauty, Shannon and Wyner. However, as mentioned above, we do not know if such large hyperbolic codes can be realized by negative curves of small genus on surfaces. In the course of proving these results, we will show the following:

\begin{thm}
\label{thm:second}
Let $\mathcal{F}$ be a set of at least two distinct irreducible curves $C$ on $X$ for which $C^2 < 0$. If $\mathcal{F}$ contains more than $R_{\rho(X) - 1}(\pi / 2)$ elements, there are two elements $C_1, C_2 \in \mathcal{F}$ together with positive integers $a, b$ such that $aC_1 + bC_2$ is an effective connected nef divisor of positive self-intersection on $X$.
\end{thm}

It was shown in \cite[Lem.\ 4.1]{MSVZ} that if the set $\mathcal{F}$ in Theorem \ref{thm:second} has more than $\rho(X)^2 + \rho(X) + 1$ elements, then there is an effective nef divisor supported on the union of the elements of $\mathcal{F}$. However, using this replaces the genus bound $b_1(X) / 4$ in Theorem \ref{thm:first} by the weaker bound $b_1(X) / (2 \rho(X)^2 + 2 \rho(X) + 2)$. In particular, it is crucial for the proof of Theorem \ref{thm:first} that we reduce down to two the number of curves in $\mathcal{F}$ involved in an effective connected nef divisor with positive self-intersection, which is clearly optimal.

\medskip

\noindent
We now outline the contents of the paper and the proofs of Theorems \ref{thm:second} and \ref{thm:first}.

\medskip

In \S \ref{s:spherical} we recall the definition of spherical codes and some classical results concerning them. We define hyperbolic codes in \S \ref{s:hyperbolic}. In \S \ref{s:goneg} we state our results concerning the relation between negative curves of small genus and hyperbolic codes of angle at least $\pi/2$.

The proof of Theorem \ref{thm:second} involves the following steps. In \S \ref{s:upperhalf} we study subsets $\mathcal{D} = \{D_i\}_i$ of $\Num(X)$ for which there is a class $h \in \Num(X)$ with $h^2 > 0$ such that $D_i^2 < 0 \le D_i \cdot D_j$ and $h \cdot D_i > 0 \ge (a D_i + b D_j)^2$ for all $i \ne j$ and all integers $a, b \ge 0$. We show $\mathcal{D}$ has these properties if and only if it forms a strict hyperbolic code with angle at least $\pi/2$ in the hyperbolic space of dimension $n = \rho(X) - 1$ associated with the intersection pairing on $\Num(X)$.

Now, suppose that $\mathcal{F}$ is a set of curves as in Theorem \ref{thm:second}. The map sending $C \in \mathcal{F}$ to its class $[C]$ in $\Num(X)$ is injective, since $C_1 \cdot C_2 \ge 0 > C_1^2$ if $C_1$ and $C_2$ are distinct elements of $\mathcal{F}$. Since $\mathcal{F}$ has at least two distinct elements, we must have $n + 1 = \rho(X) \ge 2$. Suppose now that $\# \mathcal{F} \ge R_{n}(\pi/2)$ elements, so that $\mathcal{D} = \{[C]:C \in \mathcal{F}\}$ has more than $R_n(\pi/2)$ elements. Taking $h$ to be the class of an ample effective divisor, we conclude that there are two curves $C_1, C_2 \in \mathcal{F}$ and integers $a,b \ge 0$ such that $E = aC_1 + bC_2$ has $E^2 > 0$. Since the $C_i$ are irreducible, we can adjust $a$ and $b$ so that $E$ becomes an effective connected nef divisor of positive self-intersection. This will prove Theorem \ref{thm:second}.

To prove Theorem \ref{thm:first}, we now let $\mathcal{F}$ be the set of irreducible curves $C$ on $X$ with $C^2 < 0$ and $g(C) < b_1(X)/4$. Suppose $\mathcal{F}$ has more than $R_n(\pi/2)$ elements. We show in \S \ref{s:negtocodes} that this leads to a contradiction in the following way.

Theorem \ref{thm:second} implies there are $C_1, C_2 \in \mathcal{F}$ and $a,b \ge 0$ such that $E = aC_1 + bC_2$ is an effective connected nef divisor of positive self-intersection. An \'etale Lefschetz theorem (see \cite{SGA2}, \cite[\S 2]{Bost}) implies that the induced homomorphism of \'etale fundamental groups
\[
\piet(|E|, x) \to \piet(X, x)
\]
at some geometric point $x \in E$ has image of finite index in $\piet(X,x)$, where $|E|$ is the reduction of $E$. In Theorem \ref{thm:JacSurjective} of \S \ref{s:negtocodes} we use a motivic weight argument to show that the natural morphism
\[
\mathrm{Jac}(C_1^\sharp) \oplus \mathrm{Jac}(C_2^\sharp) \to \Alb(X)
\]
is surjective, where $\mathrm{Jac}(C_i^\sharp)$ is the Jacobian of the normalization $C_i^\sharp$ of $C_i$ and $\Alb(X)$ is the Albanese variety of $X$. Since $g(C_i) = g(C_i^\sharp) = \mathrm{dim}(\mathrm{Jac}(C_i^\sharp))$
this implies that
\[
\mathrm{max}(g(C_1), g(C_2)) \ge \mathrm{dim}(\Alb(X))/2 = b_1(X) / 4.
\]
This is impossible since the curves $C$ in Theorem \ref{thm:first} are assumed to have geometric genus strictly less than $b_1(X)/4$.

This reduces the proof of Theorem \ref{thm:first} to bounding $R_n(\pi/2)$ from above. We do this in \S \ref{s:upperhalf} and \S \ref{s:mnupper} using the upper half-space model of hyperbolic $n$-space to connect hyperbolic codes to spherical codes. The connection comes about from the fact that geodesic half-spaces in the upper half-space model are one side of either a vertical plane or a Euclidean sphere. The latter spheres have centers on the Euclidean space $\mR^{n-1}$ of points at infinity different from $\infty$. In \S \ref{s:mnupper} we prove an upper bound for the size of a strict hyperbolic code of angle $\pi/2$ by showing that we can assume all the half-spaces associated with elements of the code have boundaries that are Euclidean spheres, and we can place the center of the smallest such sphere at the origin in $\mR^{n-1}$. Then the rays outward from the origin to the centers associated to other spheres must intersect a unit sphere $\mathbb{S}^{n-2}$ in $\mR^{n-1}$ in a spherical code with angle at least $\mathrm{arccos}(3/4) = \phi_0$. In \S \ref{s:mnlower} we prove a lower bound on the maximum size of a hyperbolic code with angle at least $\pi/2$ by placing the above centers at a well-chosen subset of a spherical code in $\mR^{n-1}$ with angle at at least $2 \phi_0$ and by taking the radii of all the spheres around these centers to be $\sqrt{7/8}$.

\medbreak

\noindent {\bf Acknowledgements}: The first author would like to thank the IHES, IMPA, and the University of Leiden for their support during the writing of this paper.

\section{Spherical codes}
 \label{s:spherical}

In this section we recall some definitions and results concerning spherical codes. See \cite{Ericson} for further details.

\begin{dfn}
\label{def:spherecode}
A \emph{spherical code} is a subset $\mathcal{S}$ of the unit sphere $\mathbb{S}^{n-1}$ in $n$-dimensional Euclidean space $\mR^n$. If $x, y \in \mathcal{S}$, the \emph{angle} $\phi(x,y)$ between $x$ and $y$ is the unique number in the range $0 \le \phi(x,y) \le \pi$ such that $\mathrm{cos}(\phi(x,y)) = x \cdot y$. Define the \emph{angle $\phi(\mathcal{S})$ of $\mathcal{S}$} to be the infimum of $\phi(x,y)$ over all distinct $x, y \in \mathcal{S}$. For $0 < \phi \le \tau \le \pi$ define
\begin{equation}
\label{eq:Ktwodef}
K_n(\phi,\tau) = \mathrm{max}\{\# \mathcal{S}\ :\ \phi \le \phi(x,y) \le \tau \ \textrm{for all}\ x, y \in \mathcal{S}\ \textrm{with}\ x \ne y\}.
\end{equation}
We set $K_n(\phi) = K_n(\phi,\pi)$.
\end{dfn}

\begin{example}
The kissing number $K_n = K_n(\pi/3)$ is the maximum number of spheres of a given positive radius that can touch a sphere of the same radius without having overlapping interiors.
\end{example}

The following result is due to Kabatiansky and Levenshtein \cite{KL}:

\begin{thm}
\label{thm:KLbound}
Suppose $0 < \phi \le \pi/3$. If $n$ is sufficiently large, then $K_n(\phi) \le c(\phi)^n$, where
\begin{equation}
c(\phi) = \frac{1}{2^{\, 0.099}\sqrt{1-\mathrm{cos}(\phi)}}.
\end{equation}
\end{thm}

In fact, \cite{KL} shows that the same conclusion holds for $\phi$ in the range from $0$ to a number slightly larger than $\pi/3$. The following lower bound is due Chabauty, Shannon and Wyner \cite[\S 1.6]{Ericson}.

\begin{thm}
\label{thm:lowerboundk}
Suppose $0 < \phi < \pi/2$ and $1 < c < \frac{1}{\mathrm{sin}(\phi)}$. If $n$ is sufficiently large, then
\begin{equation}
K_n(\phi) \ge c^n.
\end{equation}
\end{thm}

\section{Hyperbolic codes}
\label{s:hyperbolic}

The hyperbolic variant of spherical codes developed in this section is motivated by the following observation. A point $w$ in a spherical code $W \subset \mathbb{S}^{n-1}$ determines and is determined by the geodesic half-space
\[
Z(w) = \{w' \in \mathbb{S}^{n-1}\ :\ \langle w, w' \rangle \le 0\},
\]
where $\langle \ , \ \rangle$ is the usual Euclidean inner product. One can thus reformulate spherical codes as collections of geodesic half-spaces of $\mathbb{S}^{n-1}$ whose outward normals form at least a certain angle at their intersections.

This interpretation carries over directly to hyperbolic space. One complication is that in hyperbolic space, half-spaces may not intersect and one half-space can properly contain another. To formulate a precise definition, we first recall the definition of the hyperboloid model $L^n$ and the ball model $B^n$ of hyperbolic $n$-space. See \cite[\S 3.2]{Ratcliffe} for details.
 
Let $\langle \ , \ \rangle: \mR^n \otimes \mR^n \to \mR$ be the usual Euclidean inner product with norm $\|\ \|^2: \mR^n \to \mR$. Define an inner product $I:(\mR^n \perp \mR) \oplus (\mR^n \perp \mR) \to \mR$ by
\[
I((v;u),(v';u')) = -\langle v, v'\rangle + u\cdot u'.
\]
For $q = (v;u)$ write $q^2 = I(q,q)$. The (upper) hyperboloid model of hyperbolic space is then
\[
L^n = \left\{ q = (v; u) \in \mR^n \perp \mR \ :\ q^2 = -\| v \|^2 + u^2 = 1 \ \textrm{and} \ u > 0 \right\}.
\]
The line element for $L^n$ is $ds = \sqrt{(dv)^2 - (du)^2}$.

Let $\underline{0}$ be the origin in $\mR^n$. Projection
\[
\pi : (\mR^n \perp \mR) \smallsetminus \{(\underline{0}; -1) \} \to \mR^n
\]
from the point $(\underline{0}; -1)$ identifies $L^n$ with the ball model
\[
B^n = \{v \in \mR^n \ :\ \|v\|^2 < 1\}
\]
of hyperbolic space. The ideal boundary of $B^n$ is the closed unit sphere $\partial B^n = \mathbb{S}^{n-1}$. Define $\overline{B}^n = B^n \cup \partial B^n$, which is the closed unit ball in $\mR^n$.

Suppose $w \in \mR^n \perp \mR$ is a negative vector, i.e., that $I(w,w) = w^2 < 0$. The ray determined by $w$ is
\[
r(w) = \{ t w\ :\ 0 < t \in \mR \},
\]
and the set of points
\begin{equation}\label{eq:defineit}
Y(w) = \{ \pi(q) \in B^n\ :\ q = (v;u) \in L^n \ \textrm{and} \ I(w,q) \le 0 \}
\end{equation}
is a closed geodesic half-space in $B^n$. Let
\[
W(w) = \{ \pi(q) \in B^n\ :\ q = (v;u) \in L^n \ \textrm{and} \ I(w,q) = 0 \}
\]
be the boundary of $Y(w)$ in $B^n$ and define $\partial Y(w) \subset \mathbb{S}^{n-1}$ (resp.\ $ \partial W(w) \subset \mathbb{S}^{n-1}$) to be the ideal boundary of $Y(w)$ (resp.\ $W(w)$). Then set $\overline{Y}(w) = Y(w) \cup \partial Y(w)$ and $\overline{W}(w) = W(w) \cup \partial W(w)$. The following is well-known (e.g., see \cite[\S 2.3]{Thurston}).

\begin{lemma}
\label{lem:bijection}
The map identifying the ray $r(w)$ with the closed half-space $Y(w)$ defines a bijection between the set of negative rays in $\mR^n \perp \mR$ and the set of closed geodesic half-spaces in $B^n$. The set $\overline{W}(w)$ is the intersection of the closed unit ball $\overline{B}^n$ with a Euclidean sphere or with a hyperplane of dimension $n-1$. If $n \ge 2$ then $\overline{W}(w)$ intersects $\mathbb{S}^{n-1}$ at right angles and $\overline{W}(w) \cap \mathbb{S}^{n-1} = \partial W(w)$ is a Euclidean sphere of dimension $n-2$ and positive radius.
\end{lemma}

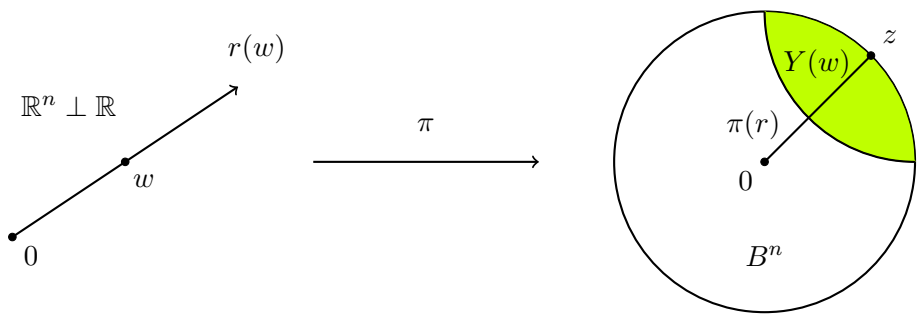
\begin{figure}\label{fig:hyperplane}
\begin{center}
\begin{tikzpicture}
\draw [thick] [->] (0,1) -- (3,3);
\draw [fill] (0,1) circle [radius=0.05];
\node at (0.25, 0.75) {$0$};
\draw [fill] (1.5, 2) circle [radius=0.05];
\node at (1.75, 1.75) {$w$};
\node at (0.75, 2.75) {$\mR^n \perp \mR$};
\node at (10, 0.75) {$B^n$};
\node at (3.25, 3.5) {$r(w)$};
\draw [black, thick] [->] (4,2) -- (7,2);
\node at (5.5, 2.5) {$\pi$};
\draw [black, thick] (10,2) circle [radius=2];
\draw [fill] (10, 2) circle [radius=0.05];
\node at (9.75, 1.75) {$0$};
\node at (9.85, 2.5) {$\pi(r)$};
\draw [fill=lime] (10, 4) arc [radius=2, start angle =180, end angle = 270];
\draw [fill=lime] (10, 4) arc [radius=2, start angle =90, end angle = 0];
\draw [color=lime, thick] (10,4) -- (12, 2);
\draw [thick] (10, 4) arc [radius=2, start angle =180, end angle = 270];
\node at (10.7, 3.3) {$Y(w)$};
\draw [thick] (10,2) -- (11.4142135624, 3.41421356237);
\draw [fill] (11.4142135624, 3.41421356237) circle [radius=0.05];
\node at (11.6642135624, 3.66421356237) {$z$};
\end{tikzpicture}
\caption{Illustration of Lemma \ref{lem:bijection}.}
\end{center}
\end{figure}

\begin{dfn}
\label{def:orho}
Suppose that $q \in \overline{W}(w)$. If $q \in W(w) \subset B^{n}$, let $n_w(q)$ be the outward unit normal to $Y(w)$ at $q$ in the tangent space $T_q\,B^n$ of $q$ in $B^n$. If $q \in \partial W(w) \subset \mathbb{S}^{n-1}$ and $n \ge 2$, let $n_w(q)$ be the outward unit normal to $\overline{Y}(w) \cap \mathbb{S}^{n-1} = \partial Y(w)$ at $q$ in the tangent space $T_q\,\mathbb{S}^{n-1}$.
\end{dfn}

We now need to understand more about the properties of the subspaces associated with a pair of negative vectors. Suppose that $w_1, w_2 \in \mR^{n}\perp \mR$ are negative vectors, and in what follows set $W_i = W(w_i)$ and $Y_i = Y(w_i)$, $i = 1,2$. Suppose that there exists a point $q \in \overline{W}_1 \cap \overline{W}_2$. It is well-known that the angle $\theta(w_1,w_2)\in [0,\pi]$ between $n_{w_1}(q)$ and $n_{w_2}(q)$ satisfies
\begin{equation}
\label{eq:thecosnew}
\mathrm{cos}(\theta(w_1,w_2)) = \frac{-I(w_1,w_2)}{\sqrt{I(w_1,w_1)\cdot I(w_2,w_2)}}.
\end{equation}
See \cite[\S 3.2]{Ratcliffe}, and note that our $I$ is the negative of the form used there. If $\overline{W}(w_1) \cap \overline{W}(w_2) = \emptyset$, define $\theta(w_1,w_2) = - \infty$. We will need the following observations.

\begin{lemma}
\label{lem:hypergeom}
Suppose $w_1$ and $w_2$ are two negative elements of $\mR^{n}\perp \mR$. The following conditions are equivalent:
\begin{enumerate}

\item[(i.)] $\theta(w_1,w_2) \ge \pi/2$;

\item[(ii.)] $I(w_1,w_2) \ge 0$ and there exists a point $q \in \overline{W}_1 \cap \overline{W}_2$;

\item[(iii.)] $I(w_1,w_2) \ge 0$ and for all $0 \le a, b \in \mR$ one has
\begin{equation}
\label{eq:abinequality}
I(a w_1 + b w_2, a w_1 + b w_2) \le 0.
\end{equation}
\end{enumerate}
Now suppose that any (and hence all) of these conditions hold and that 
there exists $h \in L^n \subset \mR^n \perp \mR$ with $I(w_i, h) > 0$ for $i = 1,2$. Then $\pi(h) \notin Y_1\cup Y_2$. Let $P \in \mathbb{S}^{n-1} = \partial B^n$ be the limit point of the geodesic ray in $B^n$ starting at $\pi(h)$ and perpendicular to $W_1$. Then $P$ is a point of $\overline{Y}_1 \smallsetminus \overline{W}_1$ that is not in $\overline{Y}_2$.
\end{lemma}

\begin{proof}
We will see that the proof reduces to checking the case $n = 2$. Since conditions (i), (ii) and (iii) are invariant under scaling, we can assume that $w_1^2 = w_2^2 = -1$. The fact that (i) implies (ii) is clear from \eqref{eq:thecosnew}. If (ii) holds, then (iii) follows from expanding $I(a w_1 + b w_2, a w_1 + b w_2)$ and using \eqref{eq:thecosnew}. Similarly, (iii) implies $|I(w_1, w_2)| \leq 1$, which means that $\overline{W}_1$ and $\overline{W}_2$ meet with angle given as in \eqref{eq:thecosnew} and hence (iii) implies (i).

We now suppose that (i), (ii) and (iii) hold and that there is an $h \in L^n$ as in the last part of the lemma. The fact that $\pi(h) \notin Y_1 \cup Y_2$ follows immediately from the definition of the $Y_i$. Intersecting with the appropriate totally geodesic $B^2$ inside $B^n$, it suffices to prove the claim for $P$ in the hyperbolic plane. Then we have the geometric arrangement shown in Figure \ref{fig:normals}. If the geodesic ray $\ell$ from $\pi(h)$ intersecting $W_1$ orthogonally were to have endpoint on $\partial \overline{Y}_2$, then it would need to also meet $\overline{W}_2$.

Let $z_1$ be the point at which $\ell$ meets $W_1$ and $z_2$ the point where $\ell$ meets $\overline{W}_2$. When $q$, $z_1$, and $z_2$ are distinct, they form a triangle in $\overline{B}^2$, possibly with an ideal vertex at $z_2$, with interior angle $\theta$ at $q$ and $\pi / 2$ at $z_1$. Therefore the triangle has angle sum greater than or equal to $\pi$, which is impossible for a triangle in $\overline{B}^2$ \cite[\S 3.5]{Ratcliffe}. In the degenerate case, $q = z_1 = z_2$, and the geodesic from $\pi(h)$ to $q$ visibly makes an angle
\[
\phi < \pi - \theta \le \pi / 2
\]
with $W_1$ at $q$, and hence cannot be orthogonal to $W_1$. Since $\pi(h)$ is not in $W_1$ and $W_1$ is totally geodesic, it is also clear that the endpoint of $\ell$ cannot be in $\overline{W}_1$. This completes the proof of Lemma \ref{lem:hypergeom}.
\begin{figure}
\begin{center}
\begin{tikzpicture}
\draw [black, thick] (0,0) circle [radius=3];
\draw [fill=blue] (-3,0) arc [radius = 3, start angle = 180, end angle = 360];
\draw [fill=lime] (2.12132034356, 2.12132034356) arc [radius = 3, start angle = 45, end angle = 225];
\draw [fill=red] (-3,0) arc [radius = 3, start angle = 180, end angle = 225];
\draw [red, fill=red] (-3,0) -- (-2.12132034356, -2.12132034356) -- (0,0);
\draw [black, thick] (3,0) -- (-3,0);
\draw [black, thick] (2.12132034356, 2.12132034356) -- (-2.12132034356, -2.12132034356);
\draw [black, ->] (0,0) -- (0,1);
\node at (0,1.25) {$w_1$};
\draw [black, ->] (0,0) -- (0.70710678118,-0.70710678118);
\node at (1, -1) {$w_2$};
\node at (-1.9, -0.7) {$Y_1 \cap Y_2$};
\node at (2, 0.5) {$\pi(h)$};
\draw [fill] (2.5,0.5) circle [radius=0.05];
\node at (-1.5, 1.5) {$Y_2$};
\node at (0, -2) {$Y_1$};
\draw [fill] (0,0) circle [radius=0.05];
\node at (-0.2, 0.25) {$q$};
\draw [thin] (0, 0.5) arc [radius = 0.5, start angle = 90, end angle = -45];
\node at (0.65, 0.3) {$\theta$};
\node at (2.2, -0.3) {$W_1$};
\node at (1.3, 1.8) {$W_2$};
\end{tikzpicture}
\caption{Geometric picture for Lemma \ref{lem:hypergeom}.}\label{fig:normals}
\end{center}
\end{figure}
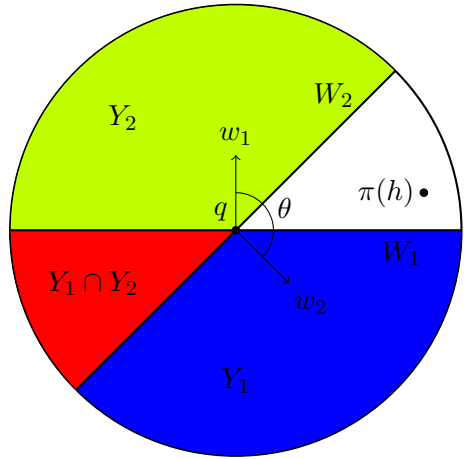
\end{proof}

\begin{dfn} 
\label{def:hypercodedef}
A \emph{hyperbolic code} is a collection $\mathcal{S}$ of negative vectors $w \in \mathbb{R}^n \perp \mathbb{R} $. We say that $\mathcal{S}$ is \emph{strict} if the union over all $w \in \mathcal{S}$ of the half-spaces $Y(w)$ is not all of $B^n$. Define $\theta(\mathcal{S}) \in \{-\infty\} \cup [0,\pi]$ to be the greatest lower bound over all pairs $w_1,w_2$ of distinct elements of $\mathcal{S}$ of the angle $\theta(w_1, w_2)$ defined above.
\end{dfn}

\begin{dfn}
\label{def:hyperdefinition}
Let $\theta$ be an angle in the range $0 < \theta \le \pi$. The \emph{hyperbolic kissing number} (resp.\ \emph{strict hyperbolic kissing number}) $\overline{R}_n(\theta)$ (resp.\ $R_n(\theta)$) in $\mZ \cup \{\infty\}$ is the supremum of $\# \mathcal{S}$ over all hyperbolic codes $\mathcal{S}$ (resp.\ strict hyperbolic codes $\mathcal{S}$) for which $\theta(\mathcal{S}) \ge \theta$.
\end{dfn}

\section{Negative curves, hyperbolic codes, and Theorem \ref{thm:first}.}
\label{s:goneg}

As in the introduction, let $X$ be an irreducible smooth projective surface over a field $k$. The group $\Num(X)$ is torsion free and finitely generated. Let
\[
\Num(X)_{\mR} = \mR \otimes_{\mZ} \Num(X).
\]
The Hodge index theorem implies that the intersection pairing on $\Num(X)$ extends to a pairing
\[
I: \Num(X)_{\mR} \times \Num(X)_{\mR} \to \mR
\]
with signature $(1,n)$, where $\mathrm{dim}(\Num(X)_{\mR}) = n + 1$.

\begin{dfn}
Let $L(X)$ be the hyperboloid model of hyperbolic $n$-space associated with the choice of isometry carrying $I$ to the standard signature $(1,n)$ pairing on $ \mR^n \perp \mR$.
\end{dfn}

\begin{dfn}
Let $\mathcal{T}(X)$ be the set of all irreducible curves $C$ on $X$ for which $C^2 = I(C,C) < 0$ and $g(C) < b_1(X)/4$. For $C \in \mathcal{T}(X)$, let $[C]$ be the class of $C$ in $\Num(X)_{\mR}$ and define $\mathcal{S}(X) = \{[ C ] : C \in \mathcal{T}\}$.
\end{dfn}

The following is a well-known consequence of the fact that distinct curves have non-negative intersection.

\begin{lemma}
\label{lem:Tbijec}
The map $\mathcal{T}(X) \to \mathcal{S}(X)$ sending $C$ to $[C]$ is a bijection.
\end{lemma}

We prove the following theorem, which is the key technical result connecting hyperbolic codes to negative curves, in \S \ref{s:negtocodes}.

\begin{thm}
\label{thm:TtoS}
The set $\mathcal{S}(X)$ is a strict hyperbolic code in $\Num(X)_{\mathbb{R}}$ of angle at least $\pi/2$.
\end{thm}

Recall that $\rho(X)$ is the rank of $\Num(X)$, i.e., the real dimension of $\Num(X)_{\mR}$. We then have the following conclusion.

\begin{cor}
\label{cor:boundnumbers}
The number of elements of $\mathcal{T}(X)$ is bounded above by the strict hyperbolic kissing number $R_n(\pi/2)$, where $n = \rho(X) - 1$.
\end{cor}
 
Recall that $K_{n-1}(\theta)$ is the kissing number associated to the angle $\theta$ in the Euclidean space $\mR^{n-1}$. The first statement of the following Theorem will be proved in \S \ref{s:upperhalf}. The second statement will be proved in \S \ref{s:mnupper} and \S \ref{s:mnlower}.

\begin{thm}
\label{thm:hyperbounds}
For every $n \ge 2$, one has
\begin{equation}
\label{eq:strictandnot}
R_n(\pi/2) \le \overline{R}_n(\pi/2) \le 2 R_n(\pi/2).
\end{equation}
Let $0 < \phi < \tau \le \pi$ be any choice of constants such that $\sqrt{2}\, \mathrm{sin}(\phi/2) = \mathrm{sin}(\tau/2)$. Then
\begin{equation}
\label{eq:upperandlowers}
\mathrm{max}\left\{ \left\lfloor \frac{K_{n-1}(2\phi_0)}{2}\right\rfloor,K_{n-1}(\phi,\tau) \right\} \le R_n(\pi/2) \le K_{n-1}(\phi_0) + 2,
\end{equation}
where $\phi_0 = \mathrm{arccos}(3/4)$.
\end{thm}

From the results about spherical kissing numbers quoted in \S \ref{s:spherical} we now have the following conclusion.

\begin{cor}
\label{cor:numeric}
One has 
\[
2^{\, 0.011\, n}(1+ o(n)) \le R_n(\pi/2) \le 2^{\, 0.901\, n}(1 + o(n)),
\]
where $o(n) \to 0$ as $n \to \infty$.
\end{cor}

\begin{proof}[Proof of Theorem \ref{thm:first}]
Combine Corollary \ref{cor:boundnumbers}, Theorem \ref{thm:hyperbounds}, and Corollary \ref{cor:numeric}.
\end{proof}

\section{Hyperbolic codes and the upper half-space model}
\label{s:upperhalf}

We now use the upper half-space model
\[
H^n =\{(z_1,\ldots,z_n)\ :\ z_i \in \mR, z_n > 0\} \subset \mR^n.
\]
of hyperbolic space to give another description of hyperbolic codes.

Recall from \cite[\S 4.4]{Ratcliffe} that there is an isometry $f : L^n \to H^n$ defined in the following way. For $x \in L^n \subset \mR^n \perp \mR = \mR^{n+1}$, let $y$ be the point $(y_1,\ldots,y_{n+1})$ on the unit $(n+1)$-sphere $\mathbb{S}^{n}$ in $\mR^{n+1}$ that is on the ray from the origin in $\mR^{n+1}$ to $x$. Then $f(x)$ is the unique point $z = (z_1,\ldots,z_n) \in H^n$ such that $(1,z_1,\dots,z_n) \in \mR^{n+1}$ lies on the ray outward from $((-1,0,0,\ldots,0);0) \in \mR^{n+1}$ through $y$.

Consider the one-point compactification
\[
\partial H^n = \{\infty\} \cup \{(z_1,\ldots,z_{n-1},0)\ :\ z_i \in \mR\}
\]
of $\mR^{n-1} = \{(z_1,\ldots,z_{n-1},0)\ :\ z_i \in \mR\}$. Then $\partial H^n$ is homeomorphic to $\mathbb{S}^{n-1}$, and the above construction identifies $\partial H^n$ with the boundary of $H^n$. Geodesics in $H^n$ are the intersection of $H^n$ with either circles or vertical lines in $\mR^n$ that intersect $\partial H^n \smallsetminus \{\infty\} = \mR^{n-1}$ orthogonally. Geodesic hypersurfaces in $H^n$ are the intersection of $H^n$ with either
\begin{itemize}

\item[(i)] vertical planes in $\mR^n$ (i.e., planes intersecting $\partial H^n \smallsetminus \{\infty\} = \mR^{n-1}$ orthogonally), or

\item[(ii)] Euclidean spheres with center on $\partial H^n \smallsetminus \{\infty\}$ (which then intersect $\partial H^n \smallsetminus \{\infty\}$ everywhere orthogonally).

\end{itemize}
Geodesic half-spaces are then formed by the set of all points of $H^n$ that lie either on one chosen side of a geodesic hypersurface or on the hypersurface itself. Define $\overline{H}^n = H^n \cup \partial H^n$.

\begin{dfn}
\label{def:upperhalf}
In \eqref{eq:defineit} to each negative vector $w \in \mR^n \perp \mR$ we defined a geodesic half-space $Y(w)$ in the open ball model $B^n$ of hyperbolic space with boundary $W(w)$, a geodesic hypersurface. Let $Y'(w)$ be the corresponding geodesic half-space in $H^n$ with boundary $W'(w)$. Similarly, let $\partial Y'(w) \subset \partial H^n$ (resp.\ $ \partial W'(w) \subset \partial H^n$) be the ideal boundary of $Y'(w)$ (resp.\ $W'(w)$). Finally, set $\overline{Y'}(w) = Y'(w) \cup \partial Y'(w)$ and $\overline{W'}(w) = W'(w) \cup \partial W'(w)$.

If $W'(w)$ lies in a vertical plane we will say that the center $z(w)$ of $Y'(w)$ is the point $\infty$ of $\partial H^n$ and that the Euclidean radius of $\partial Y'(w)$ is $\infty$. Otherwise, $W'(w)$ is the intersection of $H^n$ with a Euclidean sphere of some positive radius $d(w)$ centered at a point $z(w) \in \mR^{n-1} = \partial H^n \smallsetminus \{\infty\}$. If $z(w) \ne \infty$ and $z(w) \in \overline{Y'}(w)$, then $Y'(w)$ is the intersection of $H^n$ with the closed Euclidean ball of radius $d(w)$ about $z(w)$. Otherwise, $Y'(w)$ is the intersection of $H^n$ with the complement of the interior of this ball.
\end{dfn}

We now reformulate the condition that $\theta(w_1,w_2) \ge \pi/2$ in Lemma \ref{lem:hypergeom} using the upper half-space model. To simplify notation in what follows, given a negative vector $w_i \in \mR^n \perp \mR$ we let $Y_i = Y(w_i)$ and similarly for the other notation from Definition \ref{def:upperhalf}.

\begin{lemma}
\label{lem:basictwo}
Suppose $w_1$ and $w_2$ are two negative elements of $\mR^n \perp \mR$ such that neither $W'_1$ nor $W'_2$ lie in a vertical plane. For $z_1, z_2, d_1, d_2$ as in Definition \ref{def:upperhalf}, let $|z_1 - z_2|$ be the Euclidean distance between $z_1$ and $z_2$ in $\mR^{n-1}$. Define $\delta_i =1$ if $z_i \in \overline{Y}'_i$, and set $\delta_i = -1$ otherwise. Then $\theta(w_1,w_2) \ge \pi/2$ if and only if and only if
\begin{align}
\sqrt{d_1^2 + d_2^2} \le &|z_1 - z_2 | \le d_1 + d_2\ &\textrm{when}&\quad \delta_1 \delta_2 = 1 \label{eq:dineq1} \\
|d_1 - d_2 | \le &|z_1 - z_2| \le \sqrt{d_1^2 + d_2^2}\ &\textrm{when}&\quad \delta_1 \delta_2 = -1. \label{eq:dineq21}
\end{align}
Finally, if $\theta(w_1,w_2) \ge \pi/2$ and $I(h,w_1), I(h,w_2) > 0$ for some $h \in L^n \subset \mR^n \perp \mR$, then $z_1 \ne z_2$.
\end{lemma}

\begin{proof}
For $i = 1,2$, the half-space $\overline{Y}'_i$ is either $\overline{H}^n \cap \overline{B}(z_i,d_i)$ (when $\delta_i = 1$) or the complement in $\overline{H}^n$ of the interior of $\overline{B}(z_i,d_i)$ (when $\delta(w_i) = -1$).

Suppose first that $\theta(w_1,w_2) \ge \pi/2$, so that there is a point $q \in \overline{W}_1 \cap \overline{W}_2$. If $\delta_1 \delta_2 = 1$, the angle between the rays from $q$ to $z_1$ and from $q$ to $z_2$ is at least $\theta(w_1,w_2) \ge \pi/2$. Therefore $|z_1 - z_2| \ge \sqrt{d_1^2 + d_2^2}$. In this case, the existence of a point in $\overline{W}'_1 \cap \overline{W}'_2$ implies that $|z_1 - z_2| \le d_1 + d_2$. This proves \eqref{eq:dineq1}.

If $\delta_1 \delta_2 = -1$, the angle between the rays from $q$ to $w_1$ and from $q$ to $w_2$ is at most $\pi/2$, rather than being at least $\pi/2$. This leads to $|z_1 - z_2| \le \sqrt{d_1^2 + d_2^2}$. Since $\overline{W}'_1$ and $\overline{W'}_2$ must intersect, we see that $d_1 + d_2 \ge |z_1 - z_2| \ge |d_1 - d_2|$. Note that $|z_1 - z_2| \le \sqrt{d_1^2 + d_2^2}$ already implies $d_1 + d_2 \ge |z_1 - z_2|$. This gives \eqref{eq:dineq21}.

For the converse, one reverses the above reasoning to show that \eqref{eq:dineq1} and \eqref{eq:dineq21} imply that $\theta(w_1,w_2) \ge \pi/2$. If $z_1 = z_2$, then $W'_1$ and $W'_2$ are the the intersection of $H^n$ with concentric spheres. Hence if $\theta(w_1,w_2) \ge \pi/2$, we would have $\overline{W}'_1 = \overline{W}'_2$ and $\theta(w_1,w_2) = \pi$. However, then $\overline{Y}_1 \cup \overline{Y}_2 = \overline{B}^n$, so there could be no $h \in L^n $ with $I(h,w_1), I(h,w_2) > 0$.
\end{proof}

We now have the following, which is one of the main technical results in this paper.

\begin{thm}
\label{thm:mainbound}
For all integers $m, n \ge 1$ , the following are equivalent:
\begin{enumerate}

\item[(1)] There are elements $w_0,\cdots,w_m \in \mR^n \perp \mR$ such that for some $h \in \mR^n \perp \mR$ one has

\medbreak

\begin{itemize}
\item[(a)] $I(h,h) > 0 > I(w_i,w_i)$,
\item[(b)] $I(h,w_i) > 0$,
\item[(c)] $I(w_i, w_j) \ge 0$, and
\item[(d)] $I(aw_i+bw_j,aw_i+bw_j) \le 0$ for all distinct $0 \le i, j \le m$ and all positive $a, b \in \mR$.
\end{itemize}

\medbreak

\item[(2)] The subset $\{w_0,\ldots,w_m\} \subset \mR^n \perp \mR$ is a strict hyperbolic code having $m+1$ elements and angle at least $\pi/2$.

\medbreak

\item[(3)] After replacing the $m+1$ element subset $\{w_0,\ldots,w_m\} \subset \mR^n \perp \mR$ by their image under an isometry, the set $\{\overline{Y}'_0, \ldots, \overline{Y}'_m\}$ of half-spaces in $H^n$ has the following description.

\medskip

\noindent
The ideal boundary of each $\overline{Y}'_i$ is a sphere centered at a point $z_i \in \mR^{n-1}$ of some radius $d_i > 0$. When $i = 0$, the point $z_0$ is the origin $\underline{0}$ of $\mR^{n-1}$, and $\overline{Y}'_0$ is the exterior in $\overline{H}^n$ of the open ball of radius $d_0$.

\medskip

\noindent
If $1 \le i \le m$, $\overline{Y}'_i$ is the intersection of $\overline{H}^n$ with the closed ball of radius $d_i$ around $z_i$. Finally, the following inequalities hold:

\medbreak

\begin{enumerate}

\item[(a)] $ |z_j|^2 > \mathrm{max}(0, d_j^2-d_0^2)$ if $1 \le j \le m$,

\item[(b)] $|d_0 - d_i| \le |z_i| \le \sqrt{d_0^2 + d_i^2}$, and

\item[(c)] $ \sqrt{d_i^2 + d_j^2} \le |z_i - z_j| \le d_i+ d_j$ if $1 \le i < j \le m$

\medbreak

\end{enumerate}
where $|z - z' |$ is the Euclidean distance between points $z, z' \in \mR$.
\end{enumerate}

\medskip

\noindent
Lastly, if $\{z_1,\ldots,z_m\}$ is any set of $m$ distinct points in $\mR^{n-1}$ for which there are positive constants $d_1, \ldots, d_m > 0$ such that condition (c) of part (3) holds, then there exist $h, w_1,\ldots, w_m \in \mR^n \perp \mR$ for which the statements in condition (1) hold for $1 \le i, j \le m$.
\end{thm}

\begin{proof}
Lemma \ref{lem:hypergeom} shows the equivalence of (1) and (2). Indeed, if $h \in \mR^n \perp \mR$ and $I(h,h) > 0$, we can replace $h$ by $h/\sqrt{I(h,h)}$ to make $h$ an element of $L^n$.

To show that (1) implies (3), let $h$ and $w_0,\ldots,w_m$ be as in (1), where as above we can assume $h \in L^n$. Let $\tilde{P}$ be the limit point on $\partial H^n$ of the geodesic ray in $H^n$ that starts at $f(h)$ and is perpendicular to the geodesic hypersurface $W'_0$, where $f : L^n \to H^n$ is the above isometry. This geodesic ray is part of a geodesic line with another limit point $\tilde{P}'$ on $\partial H^n$. Applying an isometry, we can assume that $\tilde{P} = \infty$ and $\tilde{P}'$ is the origin $\overline{0}$ of $\mR^{n-1} \subset \partial H^n$.
 
Translating the final statement of Lemma \ref{lem:hypergeom} to the upper half plane model, $\tilde{P}$ is a point of $\overline{Y}'_0 \smallsetminus \overline{W}'_0$ that does not lie on $\overline{Y}'_j$ for any $j > 0$. Consider the points $z_0,\ldots,z_m$ associated with the $w_i$ by Definition \ref{def:upperhalf}. It is clear from our assumptions that each $z_j$ lies in $\mR^{n-1} \subset \partial{H}^n$. Recall that $W'_j$ is the intersection of $H^n$ with a Euclidean sphere of radius $d_j > 0$ and center $z_j$.

When $j = 0$, we know $\infty = P' \in \overline{Y}'_0$, so $\overline{Y}'_0$ must be the complement in $\overline{H}^n$ of the interior $B(z_0,d_0)$ of the ball at $z_0$ of radius $d_0$. Thus $\delta_0 = -1$ in the terminology of Lemma \ref{lem:basictwo}. Furthermore, the sphere $W_0$ is perpendicular to the geodesic with limit points $\infty$ and $\overline{0}$, and this geodesic contains $f(h)$, and we conclude that $z_0 = \overline{0}$. Note also that now $f(h)$ must lie in the interior $\ell^0$ of the vertical line segment of Euclidean length $d_0$ that has one endpoint at $\overline{0}$. See Figure \ref{fig:h2example}.

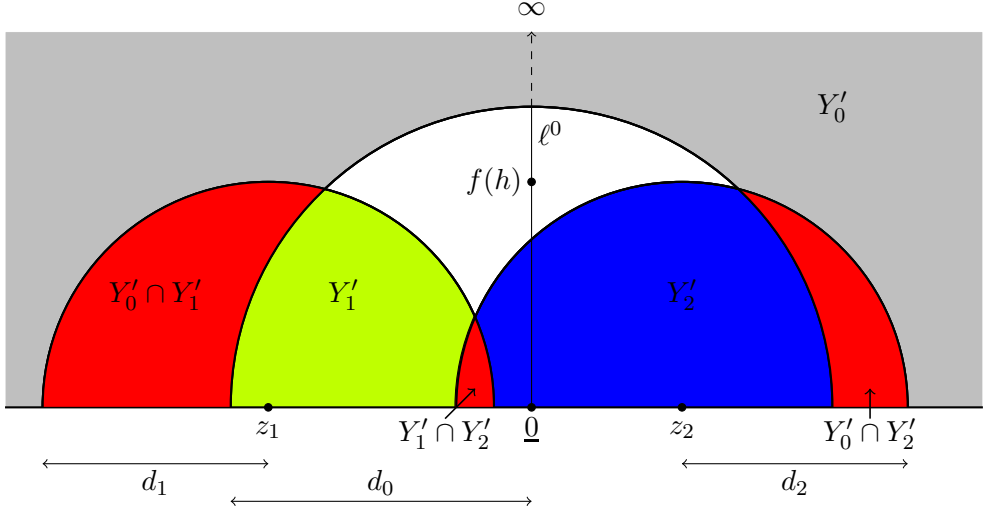
\begin{figure}
\begin{center}
\begin{tikzpicture}
\draw [white, fill=lightgray] (-7,0) -- (-7,5) -- (6,5) -- (6,0);
\draw [thick, fill=white] (4,0) arc [radius = 4, start angle = 0, end angle = 180];
\draw [thick, fill=lime] (-0.5,0) arc [radius = 3, start angle = 0, end angle = 180];
\draw [thick, fill=blue] (5,0) arc [radius = 3, start angle = 0, end angle = 180];
\draw [thick] (-0.5,0) arc [radius = 3, start angle = 0, end angle = 180];
\draw [thick] (4,0) arc [radius = 4, start angle = 0, end angle = 180];
\draw [thick] (5,0) arc [radius = 3, start angle = 0, end angle = 180];
\draw [thick] (-7,0) -- (6,0);
\draw [fill] (0,0) circle [radius=0.05];
\node at (0, -0.3) {$\underline{0}$};
\draw [fill] (-3.5,0) circle [radius=0.05];
\node at (-3.5, -0.3) {$z_1$};
\draw [fill] (2,0) circle [radius=0.05];
\node at (2, -0.3) {$z_2$};
\filldraw [fill=red] (46.56771648:4) arc (46.56771648:0:4) -- (5,0) arc (0:75.5215649:3);
\filldraw [fill=red] (136.56771648:4) arc (136.56771648:180:4) -- (-6.5,0) arc (180:75.5215649:3);
\filldraw [fill=red] (-1,0) -- (-0.5,0) arc (0:24:3) arc (156:180:3);
\draw [thick] (-0.5,0) arc [radius = 3, start angle = 0, end angle = 180];
\draw [thick] (4,0) arc [radius = 4, start angle = 0, end angle = 180];
\draw [thick] (5,0) arc [radius = 3, start angle = 0, end angle = 180];
\draw [thin, <->] (-6.5,-0.75) -- (-3.5,-0.75);
\node at (-5, -1) {$d_1$};
\draw [thin, <->] (2,-0.75) -- (5,-0.75);
\node at (3.5, -1) {$d_2$};
\node at (4,4) {$Y'_0$};
\node at (2,1.5) {$Y'_2$};
\node at (-2.5,1.5) {$Y'_1$};
\draw [fill] (0,3) circle [radius=0.05];
\node at (-0.5,3) {$f(h)$};
\draw [thin, -] (0,0) -- (0,4);
\node at (0.25, 3.65) {$\ell^0$};
\draw [thin, dashed, ->] (0,4) -- (0,5);
\node at (0,5.3) {$\infty$};
\draw [thin, <->] (0,-1.25) -- (-4,-1.25);
\node at (-2, -1) {$d_0$};
\node at (-5,1.5) {$Y'_0 \cap Y'_1$};
\node at (4.5,-0.35) {$Y'_0 \cap Y'_2$};
\draw [semithick, ->] (4.5,-0.15) -- (4.5,0.25);
\node at (-1.15,-0.35) {$Y'_1 \cap Y'_2$};
\draw [semithick, ->] (-1.15,-0.15) -- (-0.75,0.25);
\end{tikzpicture}
\caption{Arrangement of half-spaces in Theorem \ref{thm:mainbound}.}\label{fig:h2example}
\end{center}
\end{figure}

Suppose $1 \le j \le m$. Then $\infty \not \in \overline{Y}'_j$ implies that $\overline{Y}'_j = \overline{H}^n \cap \overline{B}(z_j,d_j)$, so $\delta_j = 1$ for $1 \le j \le m$. Since $f(h)$ is not contained in $\overline{Y}'_j$, we find from the fact that $f(h)$ is on $\ell^0$ that
\[
d_0^2 + |z_j|^2 > d_j^2.
\]
We now apply the criterion in Lemma \ref{lem:basictwo} to every pair $w_i, w_j$ with $0 \le i \ne j \le m$ to produce the inequalities in part (3) of Theorem \ref{thm:mainbound}.

Conversely, suppose all of the inequalities stated in part (3) of Theorem \ref{thm:mainbound} are satisfied with $z_0 = \underline{0}$ and some of $\{z_i\}_{i=1}^m \subset \mR^{n-1}$ and positive real numbers $\{d_i\}_{i = 0}^m $. Then $z_0 \ne z_i$ for $1 \le i \le m$ because $|z_i| \ne 0$ was assumed in part (3c) of Theorem \ref{thm:mainbound}. We can choose negative vectors $w_0,\ldots,w_m$ in $\mathbb{R}^n \perp \mathbb{R}$ such that $\overline{Y}'_0$ is the complement in $\overline{H}^n$ of the open unit ball of radius $d_0$ about the origin $z_0$ and $\overline{Y}'_i$ is $\overline{H}^n \cap \overline{B}(z_i,d_i)$ for $1 \le i \le m$.

The assumption that $d_0^2 + |z_i|^2 > d_i^2$ for $1 \le i \le m$ in part (3) of Theorem \ref{thm:mainbound} implies that if we choose $h \in L^n$ so that $f(h)$ lies in $B(\underline{0},1)$ and is close enough to the point that lies at distance $1$ directly above the origin, then $f(h)$ will not be in $\overline{Y}'_i$ for $1 \le i \le m$ or in $\overline{Y}'_0$. Now Lemma \ref{lem:basictwo} shows that $h,w_0,w_1,\ldots,w_m$ satisfy the conditions in part (1) of Theorem \ref{thm:mainbound}.

The final statement we must prove is that if one has only points $z_1,\ldots,z_m$ in $\mR^{n-1}$ and positive numbers $d_1,\ldots,d_m$ for which part (c) of condition (3) holds, then there are $h, w_1,\ldots,w_m \in N(X)_\mR$ for which condition (1) holds for $1 \le i, j \le m$. In this case, we choose $w_i$ so that $\overline{Y}'_i$ is $\overline{H}^n \cap \overline{B}(z_i,d_i)$ for $1 \le i \le m$. Then the vertical heights of points of each $\overline{Y}'_i$ are bounded, so we can find a point $f(h) \in H^n$ not in this union. Lemma \ref{lem:basictwo} now shows that $h,w_1,\ldots,w_m$ satisfy the conditions in part (1) of the theorem for $1 \le i, j \le m$.
\end{proof}

We now give a number of corollaries to Theorem \ref{thm:mainbound}.

\begin{cor}
\label{cor:easycor}
The strict hyperbolic kissing number $R_n(\pi/2)$ is the supremum of $m+1$ over all integers $m$ for which there exist distinct points $z_1,\cdots,z_m \in \mR^{n-1}$ and positive real constants $d_1,\cdots, d_m$ for which
\medbreak
\begin{enumerate}

\item[(a)] $\mathrm{max}\{0,d_i^2 - 1\} < |z_i|$,

\item[(b)] $|1 - d_i| \le |d_i| \le \sqrt{1 + d_i^2}$, and

\item[(c)] $ \sqrt{d_i^2 + d_j^2} \le |z_i - z_j| \le d_i + d_j$.

\end{enumerate}
\end{cor}

\begin{proof}
The corollary follows from renormalizing the $z_i$ and $d_i$ as in part (3) of Theorem \ref{thm:mainbound} by dividing each by $d_0$, $1 \le i \le m$.
\end{proof}

\begin{cor}
\label{cor:hardcor} Suppose there is a possibly nonstrict hyperbolic code in $L^n$ having $m'$ elements and angle at least $\pi/2$. If $m = \lceil m'/2 \rceil$, then there is a strict hyperbolic code having at least $m$ elements and angle at least $\pi/2$.
\end{cor}

\begin{proof}
Applying an isometry, suppose we have a nonstrict code $\mathcal{S} = \{w_1,\ldots,w_{m'}\}$ with angle at least $\pi/2$ such that, in the upper half-space model each $w_i$ gives a point $z_i \in \mR^{n-1}$ together with a positive radius $d_i$. Removing at most half the $w_i$, we can replace $m'$ by $m$ and assume that all of the constants $\delta_i$ are the same. In other words, all of the half-spaces $\overline{Y}'_i$ come from either the interiors $U_i$ of the open balls $B(z_i,d_i)$, or all of them come from the exterior of its closure of $U_i$. This means that we now satisfy the inequalities in \eqref{eq:dineq1} of Lemma \ref{lem:basictwo} for all $1 \le i,j \le m$ with $i = j$. Replacing each $\overline{Y}'_i$ by $U_i$ now produces, via the last statement of Theorem \ref{thm:mainbound}, a strict hyperbolic code with $m$ elements, since the union of all the $U_i$ cannot be all of $H^n$.
\end{proof}

\begin{cor}
\label{cor:unbound}
Suppose that there is a largest positive integer $m = m(n)$ for which the equivalent conditions (1), (2), and (3) in Theorem \ref{thm:mainbound} can be satisfied by some choice of $h, w_i, z_i$ and $d_i$ as $i$ ranges over $0 \le i \le m$. Then $m +1$ is the strict hyperbolic kissing number $R_n(\pi/2)$. The (nonstrict) hyperbolic kissing number $\overline{R}_n(\pi/2)$ satisfies
\begin{equation}
\label{eq:limits}
R_n(\pi/2) \le \overline{R}_n(\pi/2) \le 2 R_n(\pi/2).
\end{equation}
\end{cor}

\section{The upper bound on $R_n(\pi/2)$.}
\label{s:mnupper}

We begin with the following technical estimate.

\begin{lemma}
\label{lem:theta1lower}
Suppose $n \ge 2$, $z_1, z_2, z_3 \in \mR^{n-1}$, $0 < d_1 \le d_2 \le d_3$ and that
\begin{equation}
\label{eq:whatup}
\sqrt{d_i^2 + d_j^2} \le |z_i - z_j| \le d_i + d_j \quad \mathrm{for\ all}\quad i \ne j
\end{equation}
as in condition (c) of Corollary \ref{cor:easycor} (cf.\ condition 2(c) of Theorem \ref{thm:mainbound}). Then $n \ge 3$ and $z_1, z_2$ and $z_3$ are not collinear. Let $0 < \theta_1 < \pi $ be the angle of the triangle $(z_1,z_2,z_3)$ at $z_1$. Then $\theta_1 \ge \phi_0 = \mathrm{arccos}(3/4)$.
\end{lemma}

\begin{proof}
Considering the subspace spanned by $z_1, z_2$, and $z_3$, we can reduce to the case where $n \le 3$. If $z_1$, $z_2$, and $z_3$ are collinear, we can assume $n = 2$ and $z_1 < z_2 < z_3$ in $\mathbb{R}^{n-1} = \mathbb{R}$. Then \eqref{eq:whatup} leads to a contradiction. Therefore after a translation and scaling, we can assume $n = 3$, $z_1 = (0,0) = \underline{0}$ is the origin in $\mathbb{R}^2$ and $0 < d_1 \le d_2 \le d_3 = 1$. The input of $z_2, z_3, d_1$ and $d_2$ is now specified by $6$ real variables, and we want to maximize the function $\mathrm{cos}(\theta_1)$ of these variables. It is a lengthy but elementary calculus exercise to show that the maximum is obtained when $\mathrm{cos}(\theta_1) = 3/4$. We list the steps involved here and include complete details in an appendix.

Regarding $d_1$, $d_2$ and $d_3 = 1$ as fixed for the moment, let $S(d_1,d_2)$ be the set of $(z_1,z_2,z_3) = (\underline{0},z_2, z_3)$ that satisfy \eqref{eq:whatup}. One checks that $\mathrm{cos}(\theta_1)$ is a continuous function on the compact set $S(d_1, d_2)$, so that it attains its maximum at some point $(z_1,z_2,z_3) = (\underline{0},z_2, z_3)$ in $S(d_1,d_2)$. To prove the lemma it suffices to show that $\theta_1 \ge \phi_0$.

The main fact we can now use is that since $(z_1, z_2, z_3) \in S(d_1,d_2)$ maximizes $\mathrm{cos}(\theta_1)$, we cannot move $z_1$, $z_2$ and $z_3$ in $\mathbb{R}^2$ and then translate $z_1$ back to $\underline{0}$ in such a way that the inequalities \eqref{eq:whatup} still hold with the same $d_1, d_2$ and $d_3 = 1$ but with a smaller value for $\theta_1$. By considering such moves, we show in the appendix that the maximum value of $\mathrm{cos}(\theta_1)$ over all possible choices of $0 < d_1 \le d_2 \le d_3 = 1$ is attained by the example in Remark \ref{rem:hereitis} below.
\end{proof}

\begin{figure}
\begin{center}
\begin{tikzpicture}
\draw (0,0) -- (2,0) -- (1.5, 1.32287565553) -- (0,0);
\draw [fill] (0,0) circle [radius=0.05];
\draw [fill] (2,0) circle [radius=0.05];
\draw [fill] (1.5, 1.32287565553) circle [radius=0.05];
\draw (0,0) circle [radius=1];
\draw (2,0) circle [radius=1];
\draw (1.5, 1.32287565553) circle [radius=1];
\node at (-0.35,0) {$z_1$};
\node at (2.35,0) {$z_2$};
\node at (1.5, 1.67287565553) {$z_3$};
\end{tikzpicture}
\caption{Optimizing Lemma \ref{lem:theta1lower}.}\label{fig:Optimum}
\end{center}
\end{figure}
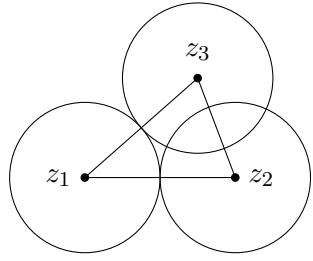

\begin{rem}
\label{rem:hereitis}
An angle of $\theta_1 = \phi_0$ can be achieved by setting $d_1= d_2 = d_3 = 1$, $n= 3$, $z_1 = (0,0) \in \mR^{n-1} = \mR^2$, $z_3 = (2,0)$ and $z_2 = (2\mathrm{cos}(\theta_0),2\mathrm{sin}(\phi_0)) = (3/2,\sqrt{7}/2)$. Then $|z_1 - z_3| = d_1 + d_3 = 2 = d_1 + d_2 = |z_1 - z_2|$ and $|z_2 - z_3| = \sqrt{d_2^2 + d_3^2} = \sqrt{2}$. See Figure \ref{fig:Optimum}.
\end{rem}

We now prove the following, which implies the upper bound in \eqref{eq:upperandlowers} of Theorem \ref{thm:hyperbounds} as well as all the bounds in \eqref{eq:strictandnot} via Corollary \ref{cor:unbound}.

\begin{cor}
\label{for:almostthere}
Suppose $z_1,\ldots,z_m \in \mR^{n-1}$ and $d_1,\ldots,d_m > 0$ satisfy condition (c) of Corollary \ref{cor:easycor}, so that $ \sqrt{d_i^2 + d_j^2} \le |z_i - z_j| \le d_i + d_j$ for all $i \ne j$. Then $m \le K_{n-1}(\phi_0) + 1$. In particular, the number $m(n)$ from Corollary \ref{cor:unbound} satisfies $m(n) \le K_{n-1}(\phi_0) + 2$.
\end{cor}

\begin{proof}
Without loss of generality, we can order $z_1,\ldots,z_m$ so that $d_1 \le d_i$ for all $1 \le i \le m$. By condition (c), the points $z_i$ are all distinct. Therefore, for $1 < i \le m$ the points
\[
\xi_i = (z_i - z_1)/|z_i - z_1|
\]
lie on the unit sphere $\mathbb{S}^{n-2}$ in $\mR^{n-1}$. Lemma \ref{lem:theta1lower} shows that for all $1 < i < j \le m$, the angle between the rays from the origin to $\xi_i$ and to $\xi_j$ must be at least $\phi_0$. Therefore $\xi_2,\ldots,\xi_m$ must form a spherical code with angular separation at least $\phi_0$, so $m-1 \le K_{n-1}(\phi_0)$. Then the number $m(n)$ from Corollary \ref{cor:unbound} is the number of points $z_0,z_1,\ldots,z_m$ for which there are $d_0,\ldots,d_m$ as in Theorem \ref{thm:mainbound}, so we conclude $m(n) \le m + 1 \le K_{n-1}(\phi_0) + 2$.
\end{proof}

\section{The lower bound on $R_n(\pi/2)$.}
\label{s:mnlower}

Let $0 < \phi < \tau \le \pi$ be any choice of constants such that $\sqrt{2}\, \mathrm{sin}(\phi/2) = \mathrm{sin}(\tau/2)$ and define $m = K_{n-1}(\phi,\tau)$. We can therefore find a spherical code $S = \{z_1,\ldots,z_m\}$ on the unit sphere $\mathbb{S}^{n-2}$ in $R^{n-1}$ such that the angular separation $\phi(z_i,z_j)$ between the rays from the origin to $z_i$ and to $z_j$ satisfies $\phi \le \phi(z_i,z_j) \le \tau$ for all $i \neq j$. Therefore,
\[
4 \mathrm{sin}^2(\phi/2) = 2 - 2 \mathrm{cos}(\phi) \le |z_i - z_j|^2 = 2 - 2 \mathrm{cos}(\phi(z_i,z_j)) \le 4 \mathrm{sin}^2(\tau/2).
\]
It follows that if we let $d_k = \sqrt{2}\, \mathrm{sin}(\phi/2) = \mathrm{sin}(\tau/2)$ for all $k = 1,\ldots,m$, then
\[
\sqrt{d_i^2 + d_j^2} \le |z_i - z_j| \le d_i + d_j
\]
for all $i \ne j$, as in condition (c) of Corollary \ref{cor:easycor}. Theorem \ref{thm:mainbound} now says that there are $h, w_1,\ldots,w_m \in N(X)_\mR$ for which the statements in condition (1) of Theorem \ref{thm:mainbound} hold for $1 \le i, j \le m$. Part (2) of Theorem \ref{thm:mainbound} now says $\{w_1,\ldots,w_m\}$ is a strict hyperbolic code with angle at least $\pi/2$. Therefore Definition \ref{def:hyperdefinition} gives that
\[
m = K_{n-1}(\phi,\tau)\le R_n(\pi/2).
\]
This is the first part of the lower bound \eqref{eq:upperandlowers} in Theorem \ref{thm:hyperbounds}.

To show the other lower bound in \eqref{eq:upperandlowers} of Theorem \ref{thm:hyperbounds}, it will suffice to show that when $\phi_0 = \mathrm{arccos}(3/4)$, we have
\[
K_{n-1}(2 \phi_0)/2 \le K_{n-1}(\phi,\tau)
\]
for some $\phi$ and $\tau$ as above. Let $\phi = 2\phi_0 = 1.445...$ and $\tau = \pi - \phi_0 = 2.418...$, so $0 < \phi < \tau \le \pi$. We then have
\begin{align*}
2 \mathrm{sin}^2(\phi/2) &= 2 \mathrm{sin}^2(\phi_0) = 2 (1 - \mathrm{cos}^2(\phi_0)) = 2 (1 - 9/16) = 7/8 \\
\mathrm{sin}^2(\tau/2) &= \mathrm{sin}^2(\pi/2 - \phi_0/2) = \mathrm{cos}^2(\phi_0/2) = \frac{ \mathrm{cos}(\phi_0) + 1}{2} = \frac{3/4 + 1}{2} = 7/8.
\end{align*}
Thus $\sqrt{2}\, \mathrm{sin}(\phi/2) = \mathrm{sin}(\tau/2)$ since both of these numbers are positive.

Recall that if $z$ and $w$ are points on the unit sphere $\mathbb{S}^{n-1}$, $\phi(z,w)$ is the angle between the rays $\tilde{z}$ and $\tilde{w}$ from the origin to $z$ and to $w$, respectively. By the definition of $\ell = K_{n-1}(2\phi_0)$, we can find a spherical code $S' = \{r_1,\ldots,r_\ell\}$ on $\mathbb{S}^{n-2}$ such that
\begin{equation}
\label{eq:boundbelowphi}
\phi(r_i,r_j) \ge 2 \phi_0\quad \mathrm{if}\quad i \ne j.
\end{equation}
For each $i$, consider the open cone $C(-r_i)$ of points $z \in \mathbb{S}^{n-2}$ such that $\phi(-r_i,z) < \phi_0$. If there were two distinct points $r_j$ and $r_q$ in $S' \cap C(-r_i)$, then
\[
\phi(r_j,r_q) \le \phi(-r_i,r_j) + \phi(-r_i, r_q) < 2\phi_0,
\]
which contradicts \eqref{eq:boundbelowphi}. Therefore there is at most point point of the form $r_j$ in $S' \cap C(-r_i)$, and if such an $r_j$ exists, $r_i$ is the unique point in $S' \cap C(-r_j)$. Throwing away at most half of the points in $S'$ we then arrive at a spherical code $S = \{z_1,\ldots,z_{\ell'}\}$ with $\ell' \ge \ell/2 = K_{n-1}(2\phi_0)/2$ such that $S \cap C(- z_i) = \emptyset$ for all $i$. If $j \ne i$, then the angle $\phi(z_i,z_j)$ can be at most $\pi - \phi_0$, since $z_j$ does not lie in $C(-z_i)$. We therefore have $\phi = 2 \phi_0 \le \phi(z_i,z_j) \le \pi - \phi_0 = \tau$, which shows that $K_{n-1}(2\phi_0)/2 \le K_{n-1}(\phi,\tau)$. This finishes the proof of the lower bound in \eqref{eq:upperandlowers} of Theorem \ref{thm:hyperbounds}.

\section{The proofs of Theorems \ref{thm:second} and \ref{thm:TtoS}.}
\label{s:negtocodes}

Theorem \ref{thm:second} is equivalent to the following result.

\begin{thm}
\label{thm:irredcurvesgood}
Suppose $\mathcal{F}$ is a set of irreducible curves $C$ on $X$ such that $C^2 < 0$ and there is no connected effective nef divisor with positive self-intersection of the form $p C_1 + q C_2$ with $0 \le p, q \in \mZ$ and $C_1, C_2 \in \mathcal{F}$. Then the set $\{[C]: C \in \mathcal{F}\}$ is a strict hyperbolic code of angle at least $\pi/2$.
\end{thm}

\begin{proof}
Recall from Lemma \ref{lem:Tbijec} that the elements $[C]$ are all distinct in $\Num(X)$. Let $A$ be an ample effective divisor on $X$. Then $I([A],[C]) > 0$ for $C \in \mathcal{F}$, where $I$ denotes the intersection pairing. Therefore $h = [A]/\sqrt{I(A,A)}$ is an element of the hyperbolic space $L(X)$ associated with the intersection pairing on $\mR \otimes_{\mZ} \Num(X)$, and it does not lie in any of the geodesic half-spaces
\[
H([C]) = \{q \in L(X): I(q,[C]) \le 0\},
\]
hence
\[
\mathcal{T} = \bigcup_{C \in \mathcal{F}} H([C])
\]
is not all of $L(X)$.

Suppose that $\mathcal{T}$ is not a strict hyperbolic code with angle at least $\pi/2$. Then we have that $\theta([C_1],[C_2]) < \pi/2$ for some distinct elements $C_1, C_2$ of $\mathcal{F}$, and Lemma \ref{lem:hypergeom} shows that there are $0 \le a, b \in \mR$ such that
\[
I(a [C_1] + b [C_2], a [C_1] + b [C_2]) = \alpha a^2 + 2\beta ab + \gamma b^2 > 0,
\]
where
\begin{align*}
\alpha &= I([C_1],[C_1]) < 0 \\
\gamma &= I([C_2],[C_2]) < 0 \\
\beta &= I([C_1],[C_2]) \ge 0.
\end{align*}
Therefore $\beta > 0$ and $\beta^2 > \alpha \gamma$. There will be positive integers $p$ and $q$ such that
\[
0 < -\gamma/\beta < p/q < -\beta/\alpha.
\]
Then $I([C_1],p[C_1] + q [C_1]) = p \alpha + \beta q > 0$ and $I([C_2],p[C_1]+ q [C_2]) = p\beta + q \gamma > 0$. It is then clear that $p C_1 + q C_2$ is an effective connected nef divisor of positive self-intersection, contradicting the hypothesis of Theorem \ref{thm:irredcurvesgood}. This proves the theorem.
\end{proof}

\begin{rem}
The referee noticed that one can also give the following simple argument for the last part of the proof of Theorem \ref{thm:irredcurvesgood}. Our setup implies that the quadratic form $\alpha a^2 + 2\beta ab + \gamma b^2$ does not take positive values on the first quadratic in $\mathbb{R}^2$. Then any positive value can be obtained with $a,b$ of the same sign, since $\alpha, \gamma < 0$ and $\beta \ge 0$, hence the form takes no positive values at all. This implies that $\beta^2 < \alpha \gamma$, which gives $\cos \theta < 0$.
\end{rem}

As in the statement of Theorem \ref{thm:TtoS}, let $\mathcal{T}(X)$ be the set of all irreducible curves $C$ on $X$ for which $C^2 = I(C,C) < 0$ and $g(C) < b_1(X)/4$, and set
\[
\mathcal{S}(X) = \{[ C ] : C \in \mathcal{T}\} \subset \mR \otimes_{\mZ} \Num(X).
\]
We must show that $\mathcal{S}(X)$ is a strict hyperbolic code in $L(X)$ with angle at least $\pi/2$. We suppose throughout this section that this is not the case, and we will derive a contradiction.

Theorem \ref{thm:irredcurvesgood} implies there is an effective connected nef divisor of positive self-intersection on $X$ of the form $pC_1+ q C_2$ in which $C_1$ and $C_2$ are elements of $\mathcal{T}(X)$ and $0 < p, q \in \mZ$. We will prove the following result below:

\begin{thm}
\label{thm:JacSurjective}
Suppose that $E$ is a connected effective nef divisor on a smooth projective geometrically integral surface $X$ over a field $k$ with positive self-intersection . Let $E^\sharp$ be the normalization of the reduction $|E|$ of $E$. Let $J(E^\sharp)$ be the direct sum of the Jacobians of the irreducible components of $E^\sharp$. Then the natural morphism from $J(E^\sharp)$ to the Albanese variety $\Alb(X)$ of $X$ is surjective.
\end{thm}

Before giving the proof, we note how it implies Theorem \ref{thm:TtoS}. If $E = pC_1+ q C_2$ as above, we obtain a surjection
\[
J(E^\sharp) = J(C_1^\sharp) \oplus J(C_2^\sharp) \to \Alb(X).
\]
Since $\Alb(X)$ has dimension $b_1(X) / 4$ and $J(C_i^\sharp)$ has dimension the geometric genus $g(C_i)$, we see that $g(C_1) + g(C_2) \ge b_1(X) / 2$. However, we supposed that every curve $C \in \mathcal{T}(X)$ has $g(C) < b_1(X)/4$, and this contradiction proves Theorem \ref{thm:TtoS}. This also completes the proof of Theorem \ref{thm:second}.

\begin{proof}[Proof of Theorem \ref{thm:JacSurjective}]
It suffices to prove the theorem for the base change of $E$ and $X$ to an algebraic closure of $k$. We assume for the rest of the proof that $k$ is algebraically closed.

Let $f$ be the pullback morphism from the Picard variety $\mathrm{Pic}^{0,red}(X)$ of $X$ to the direct sum $\mathrm{Pic}^{0,red}(E^\sharp)$ of the Picard varieties of the irreducible components of $E^\sharp$. By duality, it will be enough to show that $\Ker(f)$ is a finite group scheme. We suppose in what follows that $\Ker(f)$ is not finite and we will derive a contradiction.

Since $\Ker(f)$ is a subgroup scheme of an abelian variety, it is an extension of an abelian variety $B$ of positive dimension by a finite group scheme. Let $\ell$ be a prime different from the characteristic of $k$. Then the $\ell$-adic Tate module $T_\ell(\Ker(f))$ is isomorphic to $T_\ell(B)$, and it is a positive rank submodule of $T_\ell(\mathrm{Pic}^{0,red}(X))$. The pullback morphism from $T_\ell(\Ker(f)) = T_\ell(B)$ to $T_\ell(\mathrm{Pic}^{0,red}(E^\sharp))$ is trivial.

We know from the \'etale Lefschetz theorem that the morphism
\[
\piet(|E|, x) \to \piet(X, x)
\]
of \'etale fundamental groups at a geometric point $x$ in the support of $|E|$ is surjective, since $E$ is connected, nef, and effective. Results of this kind go back to Grothendieck in \cite{SGA2}; see Bost \cite[\S 2]{Bost} for an excellent discussion, particularly Prop.\ 2.3. This means that
\[
\Hom(\piet(X, x), \mZ/\ell^n) \to \Hom(\piet(|E|, x), \mZ/\ell^n)
\]
is injective for all $n$. Since $k$ is algebraically closed, $\mZ/\ell^n$ is isomorphic to the group scheme $\mu_{\ell^n}$ of $(\ell^n)^{th}$ roots of unity. Hence the Kummer sequence shows that
\[
\Hom(\piet(X, x), \mZ/\ell^n) = \mathrm{Pic}(X)[\ell^n] \to \Hom(\piet(|E|,x),\mZ/\ell^n) = \mathrm{Pic}(|E|)[\ell^n]
\]
is injective for all $n$.

Taking inverse limits over $n$ we see that the pullback homomorphism
\[
T_\ell(\mathrm{Pic}(X)) \to T_\ell(\mathrm{Pic}(|E|))
\]
is injective. On the other hand, $T_\ell(B) \subseteq T_\ell(\mathrm{Pic}(X))$ maps to $0$ in $T_\ell(\mathrm{Pic}^{0,red}(E^\sharp))$, so the pullback of line bundles must induce an injection
\begin{equation}\label{eq:ladic}
\xi: T_\ell(B) \to U = \Ker \left( T_\ell(\mathrm{Pic}(|E|)) \to T_\ell(\mathrm{Pic}^{0,red}(E^\sharp) \right).
\end{equation}
We will derive our contradiction from this statement.

All of the above schemes are defined over finitely generated algebras over $\mZ$. By increasing $\ell$, if necessary, we can find a specialization of all of the above schemes over a finite field $k'$ of characteristic $p$ not equal to $\ell$ so that it will suffice to show the map $\xi$ in \eqref{eq:ladic} is not injective for $k$ an algebraic closure of $k'$.

We now analyze $U$ using the map $\pi:E^\sharp \to |E|$ coming from the fact that $E^\sharp$ is the normalization of $|E|$. We have an exact sequence of sheaves of groups in the \'etale topology of $|E|$ given by
\[
1 \to \mathbb{G}_{m,|E|} \to \pi_* \mathbb{G}_{m,E^\sharp} \to V \to 1
\]
in which $V$ has support of dimension $0$. Since $\pi$ is finite, when we take the \'etale cohomology of this sequence, we find that $U$ is a quotient of
\[
M = \lim_{\overset{\longleftarrow}{n}} H^0(k\otimes_{k'} |E|,V)[\ell^n],
\]
where $H^0(k\otimes_{k'} |E|,V)[\ell^n]$ is the $\ell^n$ torsion in the $H^0(k\otimes_{k'} |E|,V)$.

Recall that $\ell$ is prime to the residue characteristic of the finite field $k'$ over which we are working. There is a filtration of $H^0(k \otimes_{k'} |E|,V)$ by $\Gal(k/k')$-stable submodules such that each graded quotient is isomorphic to either $k^*$ or the additive group $k^+$. Therefore, if $\Phi$ is the arithmetic Frobenius of $\Gal(k/k')$, then the eigenvalues of $\Phi$ on $M$ are all equal to the order $\# k'$ of $k'$. This implies that the eigenvalues of $\Phi$ on $U$ equal $\#k'$.

On the other hand $T_\ell(B)$ is the Tate module of an abelian variety $B$ over $k'$, so the eigenvalues of $\Phi$ on $T_\ell(B)$ have absolute value the square root of $\#k'$ by the Weil conjectures. It follows from this that $\xi$ cannot be injective, since $T_\ell(B)$ has positive rank. This contradiction completes the proof.
\end{proof}

\begin{rem}
We note that one can prove Theorem \ref{thm:JacSurjective} using only the classical Lefschetz theorem when $k = \mathbb{C}$, as the maps
\[
\pi_1(E^\sharp) \to \pi_1(E) \to \pi_1(X)
\]
are surjective. Our proof adapts this idea to arbitrary characteristic.
\end{rem}

\section{Appendix: A calculus exercise}
\label{s:details}

In this appendix we complete the proof of Lemma \ref{lem:theta1lower}, whose notation we now assume. As in that proof, we begin by fixing $0 < d_1 \le d_2 \le d_3 = 1$. Let $\underline{0} = (0,0)$ be the origin in $\mathbb{R}^2 = \mathbb{R}^{n-1}$, and let $S(d_1,d_2)$ be the set of triples $(z_1,z_2,z_3) = (\underline{0},z_2, z_3)$ with $z_2, z_3 \in \mathbb{R}^2$ that satisfy \eqref{eq:whatup}. Then \eqref{eq:whatup} implies that $|z_1 - z_2| = |z_2|$, $|z_1 - z_3| = |z_3|$, and $|z_2 - z_3|$ are bounded above and below by positive constants. It follows that $S(d_1,d_2)$ is compact. The law of cosines gives
\begin{equation}
\label{eq:cosinelaw}
\mathrm{cos}(\theta_1) = \frac{|z_1 - z_2|^2 + |z_1 - z_3|^2 - |z_2 - z_3|^2}{2\cdot |z_1 - z_2| \cdot |z_1 - z_3|},
\end{equation}
where the denominator on the right is bounded away from $0$. Thus $\mathrm{cos}(\theta_1)$ is a continuous function on $S(d_1,d_2)$, so it attains its maximum. We now assume this maximum occurs at $(z_1,z_2,z_3) = (\underline{0},z_2,z_3)$. As noted in \S \ref{s:mnupper}, to prove Lemma \ref{lem:theta1lower} it will suffice to show that $\theta_1 \ge \phi_0$.

Let $\theta_2$ and $\theta_3$ be the angles at $z_2$ and at $z_3$ between the sides of the triangle with vertices
at $z_1, z_2, z_3$, respectively. Since we observed in \S \ref{s:mnupper} that \eqref{eq:whatup} implies that $z_1, z_2$, and $z_3$ are not collinear, all of $\theta_1, \theta_2$ and $\theta_3$ lie in the open interval $(0,\pi)$. Suppose $\theta_3 \ge \pi/2$. Then \eqref{eq:whatup} gives
\begin{eqnarray}
\label{eq:chain}
(d_1 + d_2)^2 &\ge& |z_1 - z_2|^2 \nonumber \\
&\ge&|z_1 - z_3|^2 + |z_2 - z_3|^2 \quad\quad (\mathrm{since}\ \theta_3 \ge \pi/2) \nonumber \\
&\ge& d_1^2 + d_3^2 + d_2^2 + d_3^2.
\end{eqnarray}
This gives
\[
2 d_1 d_2 \ge 2 d_3^2 = 2.
\]
However, $d_1 \le d_2 \le d_3 = 1$, so this is only possible if $d_1 = d_2 = d_3 = 1$ and if all of the inequalities in \eqref{eq:chain} are equalities. Hence $(z_1, z_2, z_3) = (\underline{0}, z_2,z_3)$ is a right triangle with side lengths
$|z_1 - z_2| = |z_2| = d_1 + d_2 = 2$, $|z_1 - z_3| = \sqrt{d_1^2 + d_3^2} = \sqrt{2}$, and $|z_2 - z_3| = \sqrt{d_2^2 + d_3^2} = \sqrt{2}$. This means that $\theta_1 = \pi/4 \ge \phi_0$ in this case, as claimed.

As noted in \S \ref{s:mnupper}, the main fact we can now apply is that since $(z_1, z_2, z_3) \in S(d_1,d_2)$ minimizes $\theta_1$, we cannot move $z_1$, $z_2$ and $z_3$ in $\mathbb{R}^2$ and then translate $z_1$ back to $\underline{0}$ 
in such a way that the inequalities \eqref{eq:whatup} still hold with the same $d_1, d_2$ and $d_3 = 1$ but a smaller value for $\theta_1$.

We may assume that $z_3$ is a point on the positive real axis by rotating both $z_2$ and $z_3$ around $z_1 = \underline{0}$. Since $0 < \theta_1 < \phi_0 < \pi/2$, the point $z_2$ now lies in the upper right quadrant. Let $C_1$ be the circle of radius $|z_1 - z_2| = |z_2|$ around $z_2$ in $\mathbb{R}^2$. Suppose that 
\begin{equation}
\label{eq:falsify}
|z_1 - z_3| < d_1 + d_3.
\end{equation}
Let $z'_2 = z_2$ and $z'_3 = z_3$. We now move $z_1 = \underline{0}$ to a point $z'_1$ on $C_1$ that lies in the upper left quadrant and is very close to $z_1$. The only side length which changes is then $|z_1 - z_3|$, which becomes $|z'_1 - z'_3| > |z_1 - z_3|$ since $z'_3 = z_3$ lies on the positive part of the real line. Since $|z_1 - z_3| < d_1 + d_3$, all of the inequalities in \eqref{eq:whatup} will hold with the same $d_1, d_2, d_3$ when we replace $(z_1,z_2,z_3)$ by $(z'_1, z'_2, z'_3) = (z'_1,z_2,z_3)$ if $z'_1$ is a point on $C_1$ that lies in the upper left quadrant and is close enough to $z_1$. We now show that the angle $\theta'_1$ between the sides meeting at $z'_1$ of the triangle $(z'_1,z'_2,z'_3)$ satisfies $\theta'_1 < \theta_1$. This will contradict the minimality of $\theta_1$ and show that \eqref{eq:falsify} cannot hold.

Let $\theta'_3$ be the angle between the sides of $(z'_1,z'_2,z'_3)$ meeting at $z'_3 = z_3$. Since $z'_1$ lies in the upper left quadrant and on the same side of the line between $z'_2 = z_2$ and $z'_3 = z_3$ as the origin $z_1 = (0,0)$, we have $0 < \theta'_3 < \theta_3 < \pi/2$. Thus $0 < \mathrm{sin}(\theta'_3) < \mathrm{sin}(\theta_3)$. The law of sines now gives
\begin{eqnarray}
\frac{\mathrm{sin}(\theta'_1)}{\mathrm{sin}(\theta'_3)} &=& \frac{|z'_2 - z'_3|}{|z'_1 - z'_2|}\nonumber \\
 &=& \frac{|z_2 - z_3|}{|z_1 - z_2|}\nonumber\\
 &=& \frac{\mathrm{sin}(\theta_1)}{\mathrm{sin}(\theta_3)}.
\end{eqnarray}
Since $\mathrm{sin}(\theta'_3) < \mathrm{sin}(\theta_3)$ we conclude that $\mathrm{sin}(\theta'_1) < \mathrm{sin}(\theta_1)$. Since we took $z'_1$ to be close to $z_1 = (0,0)$ on $C_1$, we can ensure that that $\theta'_1$ is close to $\theta_1$. Since $0 < \theta_1 < \phi_0 < \pi/2$, we conclude that $\theta'_1 < \theta_1$, contradicting the minimality of $\theta_1$. Thus \eqref{eq:falsify} is false, so
\begin{equation}
\label{eq:simplify1}
z_1 = \underline{0} = (0,0) \quad \mathrm{and} \quad z_3 = (d_1+ d_3, 0)
\end{equation}
after rotating $z_3$ as above so that it lies on the positive real line.

Now suppose that
\begin{equation}
\label{eq:rotatez2}
d_1^2 + d_2^2 < |z_1 - z_2|^2 < (d_1 + d_2)^2.
\end{equation}
Recall that $z_2$ is a point in the upper right quadrant, and that we have reduced to the case in which \eqref{eq:simplify1} holds. We let $z'_2 = z_2$ and $z'_3 = z_3$. Define $C_2$ to be the circle with center $z_3 = (d_1 + d_3, 0)$ and radius $d_1 + d_2$, so that $C_2$ contains $z_1 = \underline{0}$ by \eqref{eq:simplify1}. Consider points $z'_1$ very close to $z_1$ on $C_2$. The only edge distance that can change on replacing $(z_1,z_2,z_3)$ by $(z'_1,z'_2,z'_3)$ is $|z_1 - z_2|$. Since $|z'_1 - z'_2| = |z'_1 - z_2|$ will be close to $|z_1 - z_2| = |z_2|$ if $z'_1$ is close to $z_1$, we conclude from \eqref{eq:rotatez2} that all the inequalities in \eqref{eq:whatup} will hold if $(z_1,z_2,z_3)$ is replaced by $(z'_1,z'_2,z'_3) = (z'_1,z_2,z_3)$ and $z'_1$ is any point on $C_2$ sufficiently close to $z_1$.

Recall that $\theta_2$ is the angle at $z_2$ between the sides of the triangle $(z_1,z_2,z_3)$ adjoining $z_2$. Let $\theta'_2$ be the corresponding angle for the triangle $(z'_1,z'_2,z'_3)$. If $z'_1$ lies in the upper half plane and is sufficiently close to $z_1$, it is on the other side of the line between $z_2$ and $z_1 = \underline{0}$ from $z_3$. It follows that $\theta'_2 > \theta_2$ in this case. We find similarly that $\theta'_2 < \theta_2$ in case $z'_1$ is a point of $C_2$ that lies in the lower half plane and is sufficiently close to $z_1$. Thus we can in either case choose a $z'_1$ on $C_2$ arbitrarily close to $z_1$ for which 
\begin{equation}
\label{eq:inequals}
0 < \mathrm{sin}(\theta'_2) < \mathrm{sin}(\theta_2).
\end{equation}
Since $|z_1 - z_3| = d_1 + d_2 = |z'_1 - z_3|$ and $|z_2 - z_3| = |z'_2 - z'_3|$, the law of sines gives
\begin{eqnarray}
 \frac{\mathrm{sin}(\theta'_1)}{\mathrm{sin}(\theta'_2)} &=& \frac{|z'_2 - z'_3|}{|z'_1 - z'_3|}\nonumber \\
 &=& \frac{|z_2 - z_3|}{|z_1 - z_3|}\nonumber\\
 &=& \frac{\mathrm{sin}(\theta_1)}{\mathrm{sin}(\theta_2)}.
\end{eqnarray}
Now \eqref{eq:inequals} shows $\mathrm{sin}(\theta'_1) < \mathrm{sin}(\theta_1)$. Since $\theta'_1$ will be close to $\theta_1 < \phi_0 < \pi/2$ for $z'_1$ close to $z_1$, we conclude that $\theta'_1 < \theta_1$, which contradicts the minimality of $\theta_1$. Thus the hypothesis \eqref{eq:rotatez2} must be false, and so
\begin{equation}
\label{eq:z2fact}
d_1^2 + d_2^2 = |z_1 - z_2|^2 \quad \mathrm{or} \quad |z_1 - z_2|^2 = (d_1 + d_2)^2.
\end{equation}

We now apply the law of cosines, together with \eqref{eq:simplify1} and $d_2^3 + d_3^2 \le |z_2 - z_3|^2$ from \eqref{eq:whatup}. This gives
\begin{eqnarray}
\label{eq:costheta1}
\mathrm{cos}(\theta_1) &=&
\frac{|z_1 - z_2|^2 + |z_1 - z_3|^2 - |z_2 - z_3|^2}{2 \cdot |z_1 - z_2| \cdot |z_1 - z_3|}
\nonumber \\
&\le&\frac{|z_1 - z_2|^2 + (d_1 + d_3)^2 - d_2^2 - d_3^2}{2 \cdot |z_1 - z_2| \cdot (d_1 + d_3)}
\end{eqnarray}
where $d_1 \le d_2 \le d_3 = 1$.

Suppose first that $d_1^2 + d_2^2 = |z_1 - z_2|$ in \eqref{eq:z2fact}. Then \eqref{eq:costheta1} becomes
\begin{equation}
\label{eq:costheta12}
\mathrm{cos}(\theta_1) \le \frac{d_1^2 + d_2^2 + (d_1+1)^2 - d_2^2 - 1}{2 \cdot \sqrt{d_1^2 + d_2^2} \cdot (d_1 + 1)}
= \frac{d_1}{ \sqrt{d_1^2 + d_2^2}} = \frac{1}{\sqrt{1 + (d_2/d_1)^2}} \le \frac{1}{\sqrt{2}}
\end{equation}
since $0 < d_1 \le d_2 $. This forces $\theta_1 \ge \pi/4$, contradicting $\theta_1 < \phi_0 < \pi/4$.

The remaining possibility in (\ref{eq:z2fact}) is that $|z_1 - z_2| = d_1 + d_2$. Then \eqref{eq:costheta1} gives
\begin{eqnarray}
\label{eq:costheta13}
 \mathrm{cos}(\theta_1) &\le&
 \frac{(d_1+d_2)^2 + (d_1 + 1)^2 - d_2^2 - 1^2}{2 \cdot (d_1+d_2) \cdot (d_1 + 1)}\nonumber \\
 &=&\frac{(d_1 + d_2 + 1)d_1 }{(d_1+d_2) \cdot (d_1 + 1)}\nonumber \\
 &\le&(1 + \frac{1}{d_1 + d_2}) \cdot (\frac{1}{1 + 1/d_1})\nonumber \\
 &\le&\frac{3}{4} \cdot 
\end{eqnarray}
since $ 0 < d_1 \le d_2 \le d_3 = 1$. This gives $\theta_1 \ge \phi_0 = \mathrm{arccos}(3/4)$, which completes the proof of Lemma \ref{lem:theta1lower}.

\bibliography{NegativeCurves}{}

\begin{thebibliography}{10}

\bibitem{AllThem}
T.~Bauer, B.~Harbourne, A.~L. Knutsen, A.~K{\"u}ronya, S.~M{\"u}ller-Stach,
  X.~Roulleau, and T.~Szemberg.
\newblock Negative curves on algebraic surfaces.
\newblock {\em Duke Math. J.}, 162(10):1877--1894, 2013.

\bibitem{Bogomolov}
F.~A. Bogomolov.
\newblock Families of curves on a surface of general type.
\newblock {\em Dokl. Akad. Nauk SSSR}, 236(5):1041--1044, 1977.

\bibitem{Bost}
J.-B. Bost.
\newblock Potential theory and {L}efschetz theorems for arithmetic surfaces.
\newblock {\em Ann. Sci. \'Ecole Norm. Sup. (4)}, 32(2):241--312, 1999.

\bibitem{ConwaySloane}
J.~H. Conway and N.~J.~A. Sloane.
\newblock {\em Sphere packings, lattices and groups}, volume 290 of {\em
  Grundlehren der Mathematischen Wissenschaften}.
\newblock Springer-Verlag, New York, third edition, 1999.

\bibitem{Ericson}
T.~Ericson and V.~Zinoviev.
\newblock {\em Codes on {E}uclidean spheres}, volume~63 of {\em North-Holland
  Mathematical Library}.
\newblock North-Holland Publishing Co., Amsterdam, 2001.

\bibitem{SGA2}
A.~Grothendieck.
\newblock {\em Cohomologie locale des faisceaux coh\'erents et th\'eor\`emes de
  {L}efschetz locaux et globaux ({SGA} 2)}.
\newblock Documents Math\'ematiques (Paris) [Mathematical Documents (Paris)],
  4. Soci\'et\'e Math\'ematique de France, Paris, 2005.
\newblock S{\'e}minaire de G{\'e}om{\'e}trie Alg{\'e}brique du Bois Marie,
  1962, Augment{\'e} d'un expos{\'e} de Mich\`ele Raynaud, With a preface and
  edited by Yves Laszlo, Revised reprint of the 1968 French original.

\bibitem{Hartshorne}
R.~Hartshorne.
\newblock {\em Algebraic geometry}.
\newblock Springer-Verlag, 1977.
\newblock Graduate Texts in Mathematics, No. 52.

\bibitem{KL}
G.~A. Kabatiansky and V.~I. Levenshtein.
\newblock Bounds for packings on sphere and in space.
\newblock {\em Problemy Peredachi Informatsii}, 14(1):3Ð25, 1978.

\bibitem{KoziarzMaubon}
V.~Koziarz and J.~Maubon.
\newblock On the equidistribution of totally geodesic submanifolds in compact
  locally symmetric spaces and application to boundedness results for negative
  curves and exceptional divisors.
\newblock arXiv:1407.656, 2015.

\bibitem{LuMiyaoka}
S.~Lu and Y.~Miyaoka.
\newblock Bounding curves in algebraic surfaces by genus and {C}hern numbers.
\newblock {\em Math. Res. Lett.}, 2(6):663--676, 1995.

\bibitem{MollerToledo}
Martin M\"{o}ller and Domingo Toledo.
\newblock Bounded negativity of self-intersection numbers of {S}himura curves
  in {S}himura surfaces.
\newblock {\em Algebra Number Theory}, 9(4):897--912, 2015.

\bibitem{MSVZ}
S.~M{\"u}ller-Stach, E.~Viehweg, and K.~Zuo.
\newblock Relative proportionality for subvarieties of moduli spaces of {$K3$}
  and abelian surfaces.
\newblock {\em Pure Appl. Math. Q.}, 5(3, Part 2):1161--1199, 2009.

\bibitem{Ratcliffe}
J.~G. Ratcliffe.
\newblock {\em Foundations of hyperbolic manifolds}, volume 149 of {\em
  Graduate Texts in Mathematics}.
\newblock Springer-Verlag, 1994.

\bibitem{Thurston}
W.~P. Thurston.
\newblock {\em Three-dimensional geometry and topology. {V}ol. 1}, volume~35 of
  {\em Princeton Mathematical Series}.
\newblock Princeton University Press, Princeton, NJ, 1997.
\newblock Edited by Silvio Levy.

\end{thebibliography}
\bibliographystyle{plain}

\end{document}